\documentclass[11pt]{article}
\usepackage{xcolor,verbatim,subcaption}
\usepackage{amssymb,amsfonts,amsmath}
\usepackage{epsfig,booktabs, authblk}
\usepackage{graphics}
\usepackage[bookmarks]{hyperref}
\hypersetup{
	colorlinks = true,
	linkcolor = blue,
	anchorcolor = blue,
	citecolor=red,
}
\usepackage[maxbibnames=1000,backend=bibtex]{biblatex}
\graphicspath{ {./figure/} }

\newcommand{\abs}[1]{\lvert#1\rvert}
\newtheorem{rmk}{Remark}
\title{Dissolution of carbonate stones caused by $CO_2$ pollutant: Numerical modelling of erosion scenarios}
\date{}

\begin{document}
	\author[1]{Elishan C. Braun}
	\author[1]{Gabriella Bretti}
	\author[2]{Samuele Ferri}
	\author[3]{Maria L. Santarelli}
	\author[4]{Matteo Semplice}
	\affil[1]{\footnotesize Institute for Applied Mathematics \textquotedblleft Mauro Picone" (IAC) – National Research Council of Italy, Rome, Italy}
	\affil[2]{\footnotesize Department of Theoretical and Applied Sciences, University of Insubria, Varese, Italy}
	\affil[3]{\footnotesize Department of Chemical Engineering Materials Environment, Sapienza University of Rome, Rome, Italy}
	\affil[4]{\footnotesize Department of Science and High Technology, University of Insubria, Como, Italy}

	\maketitle

\begin{abstract}
    In this paper, we introduce a mathematical model of carbonate-stone erosion driven by the penetration of $CO_2$-derived acidity into the water-filled pore network. The model couples transport of the reactive aqueous species, calcite dissolution and porosity evolution, thereby representing the feedback by which dissolution increases porosity and modifies subsequent transport.  
Such model is formulated as nonlinear reaction-transport system in porous media governed by Darcy flow. We propose a numerical algorithm based on finite difference approximation that relies on level-set method at the boundaries and we show numerical tests that are in accordance with the literature in terms of the advancement of the erosion front. 
\end{abstract}
	
\section{Introduction}
Among the numerous conservative issues regarding stone artifacts, here we are interested in the chemical damage caused by the presence of pollutants in the environment.  Indeed, as reported by Environmental Protection Agency, ancient buildings and sculptures in a number of cities, have weathered more during the last 20 years than in the preceding 2000. Acid rain, caused by the lowering of pH levels below 5 (acidification) of rainfalls, is mainly due to the presence of carbon dioxide ($CO_2$), sulphur oxides ($SO_x$) and, partially, nitrogen oxides $NO_x$, that increase their concentration in the atmosphere both naturally and due to human activities. 
This can be probably addressed by the fact that since the mid-19th century, human emissions along with progressive increase in the $CO_2$ concentration in the atmosphere caused by climate changes, are leading to a gradual increase of this pollutant in the environment. In this respect, NOAA Climate.gov presented studies about the increase of $CO_2$ concentration in the atmosphere\footnote{https://www.climate.gov/news-features/understanding-climate/climate-change-atmospheric-carbon-dioxide}.

Since a large number of ancient buildings and monuments are made of marble, this material has been largely studied in the last 50 years, also because it represents a simpler system than other stones thanks to its chemical homogeneity and low porosity that allows the identification of an interaction interface between the atmosphere and porous matrix, see \cite{brimblecombe}.  Marble is mainly composed of calcite, the stable polymorph of calcium carbonate, indicated in chemistry by the formula $CaCO_3$.

\begin{figure}
    \centering
    \includegraphics[width=0.4\linewidth]{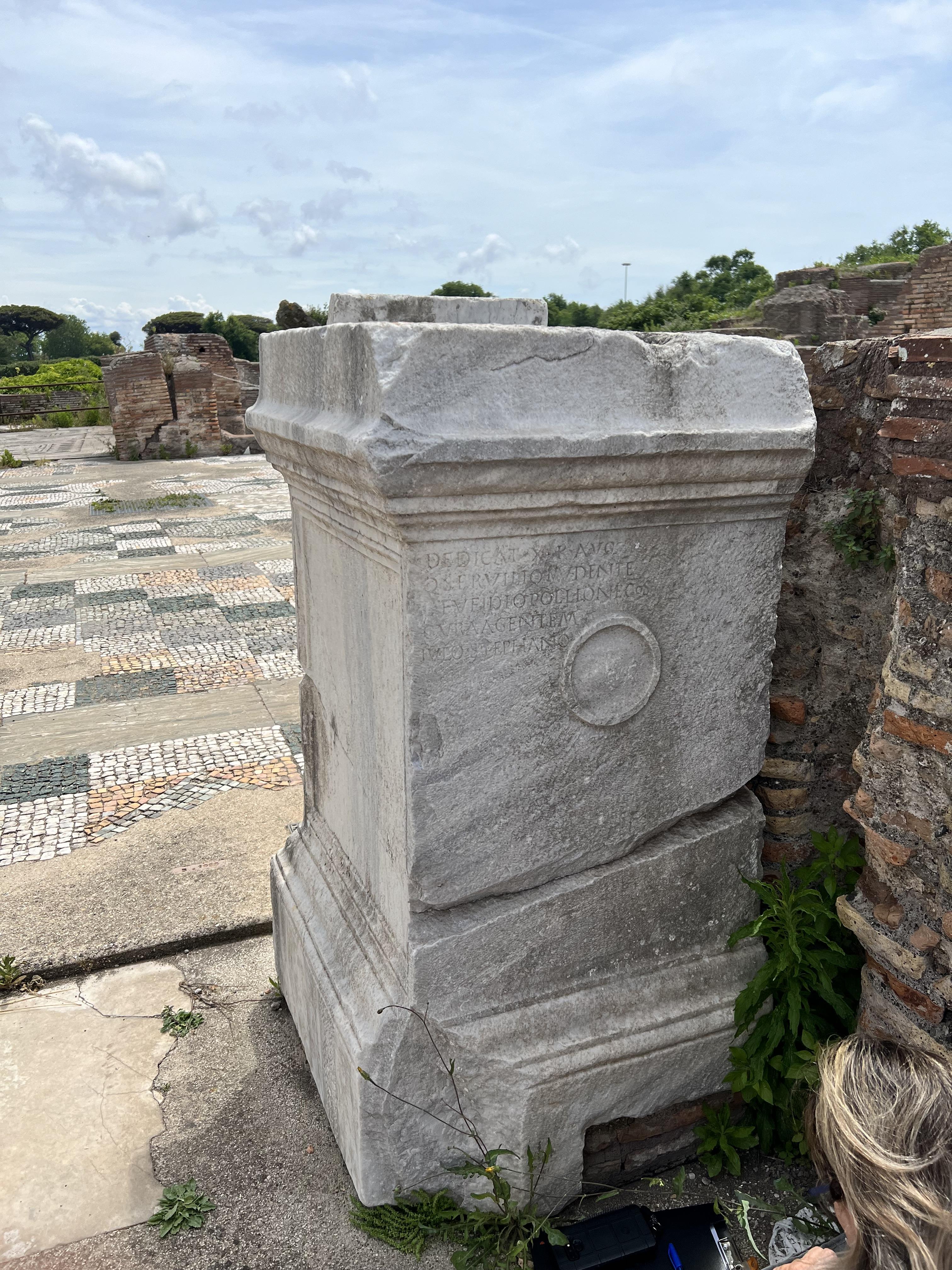}
    \caption{Marble altar from the Archaeological Park of Ostia Antica.}
    \label{fig:altare_ostia}
\end{figure}

The aim of the present work is to develop forecasting algorithm, able to correlate the pollution levels with the damage produced on carbonate materials constituting building heritage, and in particular on marble stone, see for instance Fig. \ref{fig:altare_ostia}, where the erosion is visible on the monument including the part with ancient Roman inscriptions.

\subsection{Mathematical context and related models}

The mathematical description of degradation fronts in porous materials has a long tradition, particularly in the study of concrete carbonation and stone sulfation. 
Carbon dioxide penetration in concrete is responsible for carbonation processes determining the shrinking and the deterioration of cementitious matrices, and this mechanism has been widely studied in the last decades with different approaches \cite{ashraf2016carbonation, pan} also including porosity variation \cite{Bretti}.
Among them we cite seminal works such as \cite{Aiki2013, Muntean2009} that provided a rigorous framework for modeling the penetration of reactive species and the subsequent motion of a carbonation front.  More recently, in the same framework, a moving boundary model for the carbonation in concrete has been developed and validated against data in \cite{Bretti_moving}. 

Similarly, the degradation of carbonate stones due to ${SO}_2$ aggression has been effectively described through free-boundary models \cite{Arebga2004, Clarelli2008}, focusing on the sulfation process and the formation of gypsum crusts. In this context, we mention the work  \cite{preprint26}, where predictions on the evolution of the degradation induced by marble sulfation on the real geometry of artifacts are provided.

The effects of acidic deposition on carbonate stones have been studied in laboratory and field experiments, see \cite{napap} and references therein. Indeed, acid rains corrode the stone by penetrating the pore structure and reacting with the material. In particular, calcium salts, which often serve to bind the crystal grains, are subject to dissolution by acid solutions, see \cite{brimblecombe}. 

A statistical tool (also considered by the Italian Ministry of Cultural Heritage) to evaluate damage risk, is the erosion index given by a statistical formula introduced by Lipfert \cite{lipfert}, that computes the average quantity of eroded material as a function of precipitation, deposition of $SO_2$ and concentration of $NOx$. However, this formula is not suitable for describing the chemical damage of building heritage, since it does not take into account the specificity of materials nor the time evolution of degradation. Moreover, as it provides an average estimation it neither allow to describe local microclimatic phenomena. In addition, it does not take into account the effects of the presence of $CO_2$, that is the most widespread pollutant in modern era, see the study in \cite{climate_change}. Indeed, in \cite{bonazza, bonazza2} quantitative mapping of surface recession for European carbonate buildings is provided, demonstrating that ${CO}_2$
driven chemical dissolution (the karst effect) is expected to become the dominant weathering factor in the 21st century, surpassing the historical impacts of sulfur deposition and acid rain.

In contrast to statistical indexes, important studies based on the physical-mathematical description of stone decay can be found in  \cite{Hoke2002, Hoke2004}, where a damage-based theory to quantify the chemical weathering of Carrara marble gravestones in the United Kingdom was developed. Their framework is based on the existence of an incubation period, i.e. an initial phase where the stone undergoes internal damage (forming a crust) without any physical loss of surface material and, once steady-state weathering is established, the stone surface recedes at linear speed with constant velocity that varies significantly depending on the environmental setting.

\subsection{Original contribution of our work}
Here we introduce a novel mathematical model that reproduces damage caused by erosion triggered by the penetration of pollutants transported by water causing the acidification of the porous matrix in carbonate stones. In particular, we consider the effect of carbonic acid derived from $CO_2$ emissions released in the environment, in order to directly correlate carbon dioxide pollution with physical recession in the stone. 
In this regard, data detected by sensors placed in situ, as in the Archaeological Park of Ostia Antica where we are conducting an acquisition campaign that is currently in progress, can be collected and plugged into the model. %
While in \cite{Hoke2004} environmental factors are lumped into a generic damage flux, our model uses atmospheric ${CO}_2$ as an environmental forcing for the penetration of ${CO}_2$-derived acidity into the water-filled pore network. The aqueous reactive species is transported through the porous matrix and promotes local calcite dissolution. Dissolution increases porosity, which in turn modifies the transport coefficients and facilitates further penetration of the reactive species. The present formulation is therefore intended as a macroscale effective reaction-transport model of the coupled transport-dissolution-porosity feedback associated with ${CO}_2$-driven carbonate weathering and the karst effect discussed in \cite{bonazza2, bonazza}, rather than as a complete aqueous geochemical model of the carbonate system.
 The chemistry of carbonate stone dissolution can be found in \cite{ erosion,napap} and references therein. In particular, in \cite{napap} a study conducted under the National Acid Precipitation Assessment Program  showed that for the average rain pH of 4.2, hydrogen ion deposition would contribute 17\% to the chemical erosion of marble for the flow rate conditions defined by the experimental design.

To summarize, the main contributions of the present work are the following:
\begin{itemize}
    \item We provide a chemical specialization to the dissolution process caused by carbonic acid attack, moving beyond the generic 'damage flux' used in diffusive engineering models. 
    \item Unlike traditional statistical tools such as the Lipfert index, our approach accounts for the time-dependent evolution of porosity and specific material properties, offering a physically constrained description of the degradation. 
    
    \item We introduce a robust numerical framework based on the level-set and ghost-point methods. This allows for the accurate tracking of the moving erosion front in complex two-dimensional geometries, providing a flexible tool for forecasting the physical recession of carbonate stones constituting real artifacts.
\end{itemize}

The paper is organized as follows. Section \ref{sec:model} is devoted to the description of the mathematical model for erosion of carbonate stones. In Section \ref{sec:tests} numerical tests are reported to show model's output in one and two dimensions, in particular the evolution of the erosion front.
A final section with concluding remarks and open problems ends the work.

\begin{table}
    \centering
    \resizebox{\columnwidth}{!}{%
    \begin{tabular}{|c|c|c|c|c|}
        \toprule
        Parameter  &  Description& Units & Value & Ref. \\
        \midrule
        $\mu$ & mean value water viscosity at $25^\circ$ & Poise $[g /cm\ s]$ & 8.9e-03& \cite{visc} \\
        $\rho_0$ & density of marble & $g/cm^{3}$ & 2.71 & \cite{siteConst}\\
        $D_c$ & diffusion coefficient of $HCO^-_3$ & $cm^2/s$ & 1.18e-05&\cite{sun}\\
        $K_c$ & reaction constant & $cm^3/ g\ s$ & 1.7e-03&\cite{inorganic}\\
        $r$ & recession rate of marble & $cm/s$ & 6.34e-11& \cite{brimblecombe}\\
        \bottomrule
    \end{tabular}
	}
    \caption{Parameters of the model taken from literature.}
    \label{table:stone}
\end{table}

\begin{table}
    \centering
    \begin{tabular}{|c|c|c|c|}
        \toprule
        Parameter  &  Description& Units & Value\\
        \midrule
        $\tilde{n}$ & porosity of marble&  - &   0.63\% \\
        $s_R$ & residual saturation of marble&  - & 0.227\\
        $s_S$ & maximal saturation of marble&  - &  0.884\\
        $D$ & Water diffusion rate of \eqref{B_der}&  $cm^2/s$  & 1.09e-05\\
        $D_{kP}$ & Water diffusion rate of \eqref{BkP_der}&  $cm^2/s$  & 1.09e-05\\
        $c$ & characteristic coefficient in \eqref{BkP_der} & $g / cm \ s^2$  & 34.19\\
        $K_s$ & permeability at saturation in \eqref{BkP_der} &$cm^2$ & 7.9e-09\\
        $\alpha$ & exponent in \eqref{BkP_der} & - & 0.5\\
        $\gamma$ & curvature parameter in \eqref{BkP_der} &  - & 2.0 \\
        \bottomrule
    \end{tabular}  
    \caption{Parameters of the absorption functions obtained with a fitting procedure for marble in \cite{BBDFMP}.}
    \label{table:fitting}
\end{table}

\begin{table}
    \centering
    \small
    \begin{tabular}{|c|c|c|c|c|}
        \toprule
        Parameter  &  Description & Units & Value & Reference\\
        \midrule
        $K_w$ & exchange rate with the air & $cm/s$  &  1e-02 & calibrated\\
        $K_a$ & penetration rate of carbonic acid & $cm/s$ &  1e-02 &calibrated\\
        $K_n$ & marble consumption factor& - &  1e-03  &calibrated\\
        $n_{max}$ & maximal porosity threshold &  - &   $20\%$&calibrated\\
        $\mathcal{E}$ &external humidity & - & 1.847e-03&calibrated\\
        $\mathcal{C}$ & external $CO_2$ & $g/cm^3$ & 5.5e-07 & \cite{bretticeseri}\\
        \bottomrule
    \end{tabular}
    \caption{Parameters of the model \eqref{modeleq}.}
    \label{table:kappa}
\end{table}

\section{The mathematical model describing the effects of carbonic acid in carbonate stones}\label{sec:model}
Carbon dioxide and water react producing carbonic acid:
\begin{equation}
    CO_2+ H_2 O \leftrightarrows^{k_+}_{k_-} H_2 CO_{3} 
\end{equation}
and at equilibrium one has: $k_- [H_2CO_3] = k_+ [CO_2]$. 

The process is caused by diffusion/transport mechanism regulating the movement of carbon dioxide and water inside carbonated porous material and the chemical reaction describing the related dissolution can be stated by the following equation:
\begin{equation}
    CaCO_{3(s)} + H_2 CO_{3} \longrightarrow Ca^{2+} + 2 H CO^-_{3}.
\end{equation} 
In the present simplified framework, atmospheric $CO_2$ provides the external forcing for the generation of $CO_2$-derived acidity in the aqueous pore phase. Gaseous $CO_2$ is not introduced as an independent dynamical variable. Instead, gas dissolution and hydration are treated as rapid processes related to the modeled transport and erosion time scales, and their net effect is represented through the boundary forcing and the effective aqueous concentration $c_a$. Within the pore water, $c_a$ is transported by the water flux, diffuses through the porous matrix and is consumed by reaction with the carbonate matrix. The local dissolution rate increases the porosity $n$, and the porosity-dependent transport terms in system \eqref{modeleq} provide the feedback from matrix dissolution to subsequent penetration of the reactive species. The effective reaction coefficient adopted in the present model is $K_c = 1.7 \cdot 10^{-3} cm^3 g^{-1} s^{-1}$, as reported in Table \ref{table:stone} and taken from \cite{inorganic}.

Let us consider the fraction of volume occupied by  water $\theta$ and the concentration of carbonic acid $c_a$ dissolved in water and the porosity $n$. The mathematical model derived to take into account the main aspects of phenomenon described above is the following:

\begin{equation}\label{modeleq}
	\resizebox{\columnwidth}{!}{%
    $\left\{\begin{array}{ll}
        \left\{\begin{array}{ll}
            \partial_t \theta = \nabla \cdot \left(\left(\frac{n}{\widetilde{n}}\right)^2 \nabla B\left(\frac{\theta}{n}\right)\right), &\textrm{ for } \boldsymbol{x} \in \mathbb{R}^d \textrm{ s. t. } n(\boldsymbol{x}) < n_{max}\\
            \theta = \mathcal{E}, &\textrm{ otherwise},\\
        \end{array}\right.\\[20pt]
        \left\{\begin{array}{ll}
            \partial_t(\theta c_a) = \nabla \cdot \left(c_a \left(\frac{n}{\widetilde{n}}\right)^2 \nabla B\left(\frac{\theta}{n}\right) + \theta D_c \nabla c_a \right) - K_n \rho_0\partial_t n, &\textrm{ for } \boldsymbol{x} \textrm{ s. t. } n(\boldsymbol{x}) < n_{max}\\
            c_a = \mathcal{C}, &\textrm{ otherwise},\\
        \end{array}\right.\\[20pt]
        \left\{\begin{array}{ll} 
            \partial_t n = K_c c_a(1-n), &\textrm{ for } \boldsymbol{x} \textrm{ s. t. } n(\boldsymbol{x}) < n_{max},\\
            n=n_{max}, &\textrm{ otherwise}
        \end{array}\right.
    \end{array}\right.$}
\end{equation}
where $\boldsymbol{x}\in\mathbb{R}^d$ with $d=1,2$, $\widetilde{n}$ is marble porosity (of unperturbed material), $\rho_0$ is the density of marble, $K_n$ is a factor accounting for the consumption of carbonate matrix, $D_c$ is the diffusion coefficient of $HCO^-_3$, $K_c$ is the reaction rate of $c_a$ inside the stone and $n_{max}$ is the maximum threshold for the porosity value after which the loss of porous matrix occurs, and their values are reported in Tables \ref{table:stone} and \ref{table:kappa}. 
Note that $\mathcal{E}$ represents the prescribed external humidity, whereas $\mathcal{C}$ is an effective boundary concentration representing the availability of $CO_2$-derived reactive carbon in the aqueous phase at the stone surface under the prescribed atmospheric $CO_2$ and moisture conditions. In the present formulation, $\mathcal{C}$ acts as a boundary forcing and is not the concentration of an independently modeled gaseous $CO_2$ phase. The conversion from measured atmospheric $CO_2$ concentration to $\mathcal{C}$ is treated as an effective calibrated closure, since Henry-law partitioning and gas-liquid mass transfer are not explicitly resolved. The values of  $\mathcal{E}$ and $\mathcal{C}$ used in the simulations are specified in Section \ref{sec:tests}.
We remark that these functions can be generalized to be time dependent and also include rain effect, as well as seasonal and environmental changes.
\begin{rmk}
The chemical description adopted here is intentionally reduced and should be interpreted as an effective kinetic representation of $CO_2$-driven carbonate dissolution. The adopted rate law is intended for conditions in which the pore solution remains sufficiently undersaturated with respect to calcite, so that dissolution dominates and the reduction in thermodynamic driving force near equilibrium can be neglected. Precipitation and reprecipitation of $CaCO_3$ are not represented. Calcium concentration, pH, alkalinity and the full aqueous carbonate speciation are not treated as independent dynamical variables; consequently, carbonate buffering, accumulation and transport of reaction products, and their feedback on calcite saturation are not explicitly resolved. These effects are represented, at most, implicitly through the effective reaction coefficient and its calibration. The approximation is therefore appropriate as a first-order macroscale description when product accumulation does not bring the pore water close to calcite saturation. It may become inaccurate under poorly flushed conditions, near equilibrium, or when precipitation is significant. Extending the model to include gas-liquid mass transfer, aqueous carbonate equilibria, calcium transport and saturation-dependent dissolution-precipitation kinetics is left for future work.
\end{rmk}

Since the crucial parameters for the evolution of the erosion front are $K_w$, $K_{a}$, $K_{n}$ and $n_{max}$ we need to calibrate them numerically against available data. In particular, we refer to the value reported in \cite{brimblecombe} where from the study of depth and thickness of inscriptions on tombstones in marbles a recession rates was found ranging from $3.6$-$2.8 \ mm/100y$ for urban areas and the estimated value of $1.7 \ mm/100y$ for a suburban site. Based on these values, we are assuming a value of $2 \ mm/100y$ as a reference for the recession rate $r$ in suburban sites due to the different $CO_2$ concentration levels. Note that the value of $r$ is reported in Table \ref{table:stone} in the unit measure of $cm/s$ and the recession rate is considered approximately linear as from experimental evidence \cite{brimblecombe}.

The system (\ref{modeleq}) is coupled with the following boundary conditions:
\begin{equation}\label{model:bc}
    \left\{\begin{array}{ll}
        \left(\frac{n}{\widetilde{n}}\right)^2\nabla_{n}B\left(\frac{\theta}{n}\right) &= K_w\left(\mathcal{E} - \theta\right) ,\\[10pt]
        c_a \left(\frac{n}{\widetilde{n}}\right)^2 \nabla_{n}B\left(\frac{\theta}{n}\right) + D_c \theta \nabla_{n} c &= K_a\left(\mathcal{C} - c_a\right),
    \end{array}\right.
\end{equation}

where $K_w$ and $K_a$ denote, respectively, the exchange rate of water and the exchange rate of the effective aqueous reactive species at the stone boundary.

\begin{rmk}
We clarify that $c_a$ denotes the effective aqueous concentration of $\text{CO}_2$-derived reactive carbon promoting calcite dissolution, rather than purely the undissociated carbonic acid species $\text{H}_2\text{CO}_3$. In the aqueous phase, dissolved carbon dioxide hydrates and dissociates rapidly compared to the macroscopic erosion time scales. 
Furthermore, while atmospheric $\text{CO}_2$ drives the acidification of the pore water, Henry's law partitioning ($\text{CO}_{2(g)} \rightleftharpoons \text{CO}_{2(aq)}$) and explicit gas-phase diffusion within the unsaturated pore network are not dynamically resolved as an independent two-phase system. Instead, the boundary concentration $\mathcal{C}$ represents the effective aqueous saturation concentration in equilibrium with ambient air, which can be linked to the atmospheric partial pressure $p_{\text{CO}_2}$ via an effective Henry's law relation $\mathcal{C} \approx K_H^* p_{\text{CO}_2}$. For standard atmospheric levels ($\approx 375~\text{ppm}$), this yields the reference boundary value $\mathcal{C} = 5.5 \times 10^{-7}~\text{g/cm}^3$ reported in Table~\ref{table:kappa}. The term $K_a(\mathcal{C} - c_a)$ in ~\eqref{model:bc} thus models the effective mass transfer of $\text{CO}_2$-derived acidity across the stone boundary as a calibrated thermodynamic closure.
\end{rmk}

\subsection{Definition of function \texorpdfstring{$B$}{B}}
In the first equation of system (\ref{modeleq}), also known as Richards equation \cite{richards}, an important aspect is the choice of absorption function $B$ connecting the capillary pressure with the moisture content, see \cite{bear, onofri} and references therein for a deeper discussion.

The function $B$ satisfies relation given by Darcy's law describing the flux in unsaturated media: 
\begin{equation}\label{Bdarcy}
    \nabla B\left(\frac{\theta}{n}\right)=-\frac{k\left(\frac{\theta}{n}\right)}{\mu} \left(\nabla P_c\left(\frac{\theta}{n}\right) - \rho_l g\right),
\end{equation}
where $P_c$ represents the capillary pressure, $k$ the relative permeability of the porous matrix, $\mu$ the viscosity of the fluid, $\rho_l$ is the density of the fluid (water) and $g$ the gravitational acceleration, that can be safely disregarded for specimen of small sizes.
The capillary pressure $P_c(s)$ is a decreasing function of saturation $s:=\theta/n$, while the relative permeability $k(s)$ is a non-negative increasing function of $s$ and is upper-bounded by its value at saturation, namely $s_S$ (see \cite{bear} for further details).

A possible approach to correlate capillary pressure with moisture content into the porous matrix is to introduce in Darcy's law a polynomial function with some free parameters identified through model calibration. In the following we describe the absorption functions which are introduced in \cite{bretti-belfiore, Clarelli2010}.

\subsubsection{The symmetric absorption function \texorpdfstring{$B'$}{B'}}\label{sec:NN}
The $B$ function introduced in \cite{Clarelli2010} is formulated as:
\begin{align} \label{NN}
    B(s) &= 
    \begin{cases}
        0 & s \in [0,s_R)  \\  
        -\frac{(2 D (s_R - s)^2 (s_R - 3 s_S + 2 s))}{3 (s_R - s_S)^2} & s \in [s_R,s_S] \\
        \frac{2}{3} D (s_S - s_R) & s > s_S,
    \end{cases} 
\end{align}
with  
\begin{equation}\label{B_der}
    \partial_s B(s) = B'(s) = max\left(0,-\frac{4 D (s_R - s) (s_S - s)}{(s_R - s_S)^2}\right),
\end{equation}
where $\{s_R, s_S, D\}$ is the set of model parameters to be determined: $s_R$, the minimum value for saturation ensuring the hydraulic continuity, $s_S$ the maximum value of $s$ reached at saturation, and the water diffusion rate $D:=\max_{[s_R,s_S]}{B'}$ reached at $s=\frac{s_R+s_S}{2}$ (see Table \ref{table:fitting}). We remark that $B'(s)$ is a compactly supported function for $s$ in $[s_R,s_S]$ and $D$ has the dimensions of a diffusivity. 
Then, the profiles of $B(s)$ and $B'(s)$ are obtained by letting the right and left endpoints $s_R$ and $s_S$ vary, i.e. $0 < s_R < s_S \le 1$, meaning that the realistic pore saturation is normally less than $100\%$. \\

\subsubsection{The asymmetric absorption function \texorpdfstring{$B_{kP}'$}{BkP'}}
Another possible more general formulation of $B$ function in Darcy's law consists in expressing separately the permeability function $k(s)$ and the capillary pressure $P_c(s)$, introduced in \cite{bretti-belfiore}. In particular, the permeability function $k(s)$ is a generalization of the function proposed in \cite[Chapter~4]{bear} expressed as:
\begin{equation}\label{perm_fun}
    k(s) = \left\{
    \begin{array}{ll}
        K_{s} \left(\frac{s-s_R}{s_S-s_R}\right)^\gamma, \textrm{ if } s \in [s_R, s_S],\\
        k(s_R) ,\textrm{ if } s < s_R,\\
        k(s_S), \textrm{ if } s > s_S,  
    \end{array}\right.
\end{equation}
with $K_{s}>0$ the constant of permeability at maximal saturation, possibly obtained by experimental measurements, and $\gamma>0$ a parameter to be calibrated with experimental data. Note that $k(s_R)=0$ and $k(s_S)= K_{s}$. 
For the capillary pressure, here we choose $P_c (s) =  c \frac{(s-s_S)^2}{(s-s_R)^{\alpha}} \ \textrm{ for } s \in (s_R,s_S]$, with $\alpha>0$ an exponent and $c>0$ a diffusion constant, see \cite{bretti-belfiore} for further details. Its derivative with respect to $s$ is then:
\begin{equation}\label{P'}
    P'_c (s) = -\frac{c (s-s_S) (2 s_R - 2s - \alpha s_S + \alpha s)}{(s - s_R)^{\alpha+1}}.
\end{equation}

Constants $s_R, s_S$ and $c$ are physical properties of the porous material involved and will be calibrated later on, with $s_R$ as the residual saturation and $s_S$ as the maximum saturation level. The quantity $\theta_R = s_R \cdot n$, corresponding to $s=s_R$, is the minimum value for saturation ensuring the hydraulic continuity, with $B'_{kP}(s)$ a compactly supported function in $[s_R,s_S]$.

The function is null outside the interval $[s_R,s_S]$ and has the following expression in $[s_R,s_S]$:
\begin{align} \label{BkP_der}
    B'_{kP}(s) &= max\left(0, K_s \frac{c}{\mu} \frac{(s-s_R)^{\gamma - \alpha -1}}{(s_S-s_R)^\gamma} (s-s_S) (2 s_R + s (\alpha -2) - \alpha s_S) \right),
\end{align}
with $D_{kP} = \max{B'_{kP}}$ the diffusion coefficient and  $B'_{kP}(s)$ can be integrated exactly as shown in \cite{bretti-belfiore} or numerically.

In this case the set of model parameters is given by $\{s_R, s_S, \alpha,c,K_{s},\gamma\}$, with $\gamma-\alpha-1>0$ and we indicate $D_{kP} = \max_{[s_R,s_S]}{B'_{kP}}$ the diffusion coefficient. Now, we integrate exactly $B'_{kP}(s)$ and we get the following expression:
\begin{equation} \label{BP2}
	\resizebox{\columnwidth}{!}{%
	$\widetilde{B}_{kP} (s) = \frac{K_s c (s - s_R)^{\gamma - \alpha}}{\mu (s_S - s_R)^{\gamma}} \frac{s^2 u + s v + \gamma^2 s_S (-2s_R + \alpha s_S) + \gamma z + \alpha s_S^2(\alpha^2 - 3 \alpha + 2)}{ (- \alpha^3 + 3 \alpha^2 (\gamma+1) - 3\alpha \gamma(\gamma+2) - 2\alpha + \gamma^3 + 3 \gamma^2 + 2 \gamma)},$}
\end{equation}
with 
\begin{align*} 
    u&= \alpha^3 - \alpha^2(2\gamma +3) +\alpha(5 \gamma +2) - 2\gamma(\gamma+1),\\
    v&= \resizebox{0.92\columnwidth}{!}{ $2\gamma(- s_R \alpha + 2 \alpha^2 s_S + 2 s_S - 4\alpha s_S) +  2\gamma^2 (s_R - \alpha s_S + s_S ) + 2 \alpha s_S (-\alpha^2 + 3 \alpha - 2) ,$}\\
    z&= 2 s_R^2 + \alpha s_S^2 (-2\alpha + 3) + 2s_R s_S(\alpha -2).
\end{align*}
Then, $B_{kP}(s)$ reads as:
\begin{align} \label{BkP}
    B_{kP}(s) &= 
    \begin{cases}
        0 & s \in [0,s_R] \\  
        \widetilde{B}_{kP}(s) \textrm { in \eqref{BP2}} & s \in (s_R,s_S) \\
        \widetilde{B}_{kP}(s_S)=\frac{2 K_s c \gamma (s_S - s_R)^{2-\alpha}}{\mu (- \alpha^3 + 3 \alpha^2(\gamma + 1) - 3 \alpha \gamma(\gamma +2) - 2\alpha + \gamma^3 + 3\gamma^2 + 2\gamma)} & s \in [s_S,1].
    \end{cases} 
\end{align}

\begin{rmk}\label{remark1}
    We end this paragraph with a brief discussion on the two definitions of the absorption functions, namely $B'$ and $B'_{kP}$. Due to its simplicity, $B'$ has the great advantage of depending only on 3 parameters and this simplifies its calibration against data. However, with this formulation the capillary pressure $P_c(s)$ and the permeability function $k(s)$ are not defined independently as in Darcy's law \eqref{Bdarcy}. On the contrary, $B'_{kP}$ includes a larger number of parameters providing a richer description of capillary and permeability properties of the material. Even if the two formulations provide comparable results if the only data available are imbibition curves, $B'_{kP}$ allows a more accurate reconstruction of the phenomenon in the presence of further data types, such as granulometric distribution (from MIP experiment) and permeability test at saturation, see \cite{bretti-belfiore}. 
    Moreover, we can take advantages from the separate formulation of capillary pressure and permeability in $B'_{kP}$ function, since we aim at including the gravity effect in Richards equation in future works.
\end{rmk}

\begin{figure}
    \centering
    \includegraphics[width=0.48\linewidth]{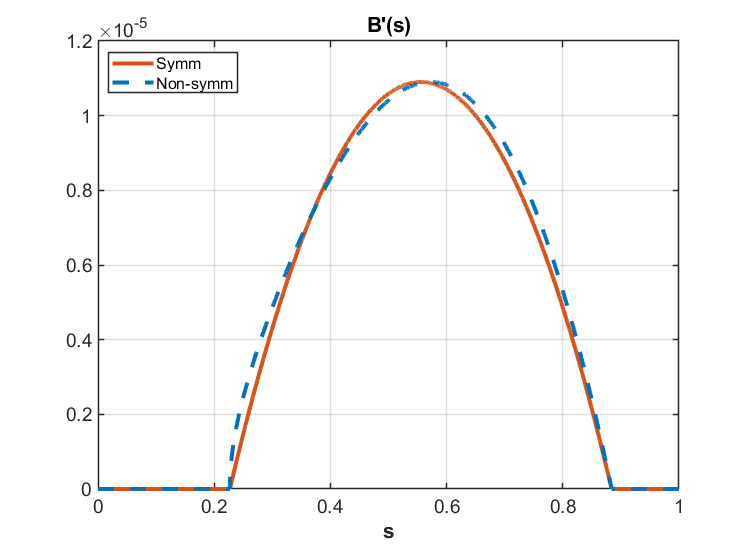}
    \includegraphics[width=0.48\linewidth]{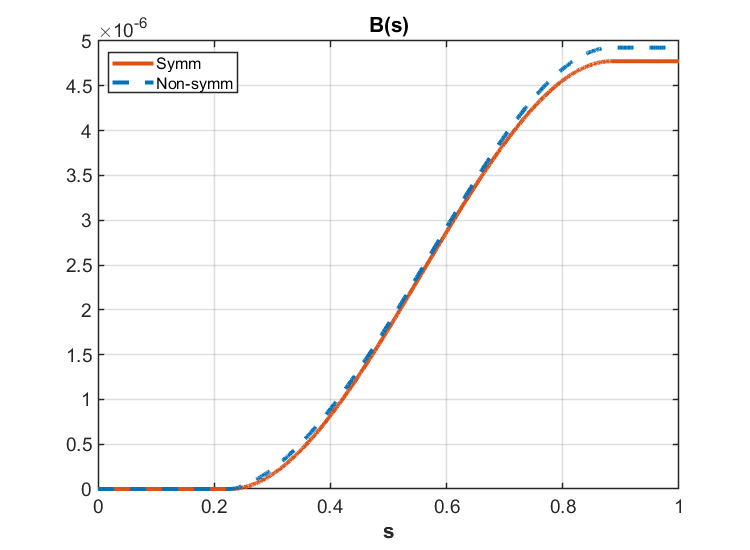}
    \caption{Profile of symmetric and non-symmetric absorption functions, $B'$ and $B'_{kP}$ (left) and of their integrands $B$ and $B_{kP}$ (right). }
    \label{fig:Bder}
\end{figure}

\subsection{Numerical approximation of the model}
In this section we describe the numerical scheme used to discretize the mathematical model \eqref{modeleq} for the two-dimensional case; the one-dimensional case is straightforward. Our scheme is a combination of finite differences at internal points of the computational domain and level-set-ghost-point method from \cite{Ghost-point,CocoFD,CocoSecondOrder} at the boundaries, since our main goal is to keep track of the moving boundary across time due to erosion.\\

Suppose $\Omega$ is the physical domain of reference (e.g. a section of a pillar). In general its boundary $\partial\Omega$ can not be aligned to a cartesian grid. Therefore we embed it in a larger cartesian domain $\overline{\Omega}\subset\mathbb{R}^2$ (a square, for instance). Here, we mesh the new domain using a uniform grid $\mathcal{G}_h=\left\{(x_i,y_j)\mid i,j=0,...,N \right\}$ with spatial step $h=x_i - x_{i-1} = y_j - y_{j-1}$, for $i,j=1, ..., N$, to obtain the computational domain $\overline{\Omega}_h$.
Now, we introduce a level-set function $\varphi:\mathbb{R}^2\rightarrow \mathbb{R}$ such that
\begin{equation}
    \Omega = \left\{\boldsymbol{x}\in\overline{\Omega} \mid\varphi(\boldsymbol{x})\leq 0\right\}\subset \overline{\Omega} \quad \text{ and }\quad \partial\Omega = \left\{\boldsymbol{x}\in\overline{\Omega} \mid \varphi(\boldsymbol{x}) = 0\right\}.
\end{equation}
In our model, the level-set function is
\begin{equation}\label{levelset}
    \varphi(\boldsymbol{x}) = n(\boldsymbol{x}) - n_{max},
\end{equation}
where $n_{max} = p \tilde{n}$ is the threshold porosity at which the loss of material occurs; it is defined as the initial porosity $\tilde{n}$ multiplied by a growth factor $p \in (1,1/\tilde{n})$ and it has to be calibrated from data.\\
We want to remark at this point that since the porosity $n(x,t)$ is a time-dependent function, the level-set function is
\begin{eqnarray}
    \varphi(\boldsymbol{x},t)=n\left(\boldsymbol{x},t\right)-n_{max}
\end{eqnarray}
with time-dependent domain $\Omega_{t}$ and boundary $\partial \Omega_{t}$, but we surpress the explicit dependence on $t$ to simplify the notation.\\
In this way we separate $\overline{\Omega}_h$ in three different subsets. Firstly, we define the set of \textit{internal points} as
\begin{equation}
    \overline{\Omega}_{h,I} := \left\{\boldsymbol{x}\in\overline{\Omega}_h \mid \varphi(\boldsymbol{x}) < 0\right\}.
\end{equation}
Then, we introduce $\overline{\Omega}_{h,B}$, the set of points in $\overline{\Omega}_h \setminus \overline{\Omega}_{h,I}$ that have a neighbor in $\overline{\Omega}_{h,I}$. We call these \textit{ghost points}.
Finally, we define the set of \textit{outside points} as
\begin{equation}
    \overline{\Omega}_{h,O} := \overline{\Omega}_h \setminus \overline{\Omega}_{h,I} \setminus \overline{\Omega}_{h,B}.
\end{equation}
For each outside point in $\overline{\Omega}_{h,O}$ we impose
\begin{equation}\label{model:ext}
    \left\{\begin{array}{ll}
        \theta = \mathcal{E} \\
        c_a = \mathcal{C} \\
        n = 1
    \end{array}\right. .
\end{equation}
In $\overline{\Omega}_{h,I}$, which contains points inside the physical space $\Omega$, we discretize the equation using a first-order 2D finite difference scheme with a 3-points wide stencil along each direction. As in \cite{Bracciale}, we use a first order approximation obtained by Taylor expansions as the simplest and consistent discretization of $\partial_x(q(x)\partial_x w(x))$ :
\begin{equation}\label{approx:eq}
    \Delta_i(q,w) := \frac{(q_i + q_{i+1})(w_{i+1}-w_i) - (q_{i-1} + q_i)(w_i-w_{i-1})}{2 \Delta x^2},
\end{equation}
so that, in the multidimensional case we can write it as:
\begin{equation}
    \begin{split}
        \nabla\cdot(q\nabla w) &= \partial_x(q\partial_x w) + \partial_y(q\partial_y w)\\
        &= \Delta_i(q,w) + \Delta_j(q,w)=: \Delta_{ij}(q,w),
    \end{split}
\end{equation}
where $\Delta_{p}(q,w)$ is the approximation \eqref{approx:eq} along the dimension $p$, namely $x$ or $y$. To simplify the notation, we write $q_{ij}^{\tau}$ to indicate the value of $q$ at point $(x_i, y_j)$ at time $\tau$.

In order to avoid strong constraint on the time step, we choose the implicit Euler method which allows large time steps. Thus the discretization schemes for the interior are:
\begin{itemize}
    \item for the porosity equation
    \begin{equation}\label{por:discr}
        \frac{n_{ij}^{\tau+1} - n_{ij}^{\tau}}{\Delta t} = K_c c_{a, ij}^{\tau+1} (1 - n_{ij}^{\tau+1}),
    \end{equation}
    
    \item for the non-linear diffusion of humidity
    \begin{equation}\label{theta:discr}
        \frac{\theta_{ij}^{\tau+1} - \theta_{ij}^{\tau}}{\Delta t}= \Delta_{ij}\left( \left(\frac{n^{\tau+1}}{\widetilde{n}}\right)^2  , B\left(\frac{\theta^{\tau+1}}{n^{\tau+1}}\right) \right),
    \end{equation}
    
    \item for the diffusion of the carbonic acid
    \begin{equation}\label{co2:discr}
        \begin{split}
            \frac{\theta_{ij}^{\tau+1}c_{a,ij}^{\tau+1} - \theta_{ij}^{\tau}c_{a,ij}^\tau}{\Delta t} =&\ \Delta_{ij}\left(c_{a,ij}^{\tau+1}\left(\frac{n_{ij}^{\tau+1}}{\tilde{n}}\right)^2, B\left(\frac{\theta_{ij}^{\tau+1}}{n_{ij}^{\tau+1}}\right)  \right) \\
            &+\Delta_{ij} \left( D\theta_{ij}^{\tau+1}, c_{a, ij}^{\tau+1} \right) - K_c c_{a,ij}^{\tau+1} (1-n^{\tau+1}_{ij}) K_n \rho_0.
        \end{split}
    \end{equation}	
\end{itemize}

However, it may happen that the finite difference stencil overflows the boundary and includes points from $\overline{\Omega}_{h,B}$, meaning that we have a system of $N_i$ equations in $N_i+N_g$ unknowns, where $N_i$ is the number of internal points, while $N_g$ is the number of external points reached by the internal stencils.
To close the system we add $N_g$ equations from the discretization of the boundary condition in each of such points.\\

Since we do not know the exact coordinates of the boundary, unless it is aligned with the grid, the strategy is to start from the ghost point and project it onto $\partial \Omega$ along the normal direction of the level-set $\varphi$, and there discretize the boundary condition through an interpolation.

We now exploit the method for the boundary condition for the humidity equation of (\ref{modeleq}). The case for the carbonic acid equation is straightforward. Fixing a ghost point $\boldsymbol{x}_G$, we project it along the normal line of the level-set to get a new point
\begin{equation}
    \boldsymbol{x}_B = \boldsymbol{x}_G + \lambda\nabla\varphi(\boldsymbol{x}_G)
\end{equation}
such that $\varphi(\boldsymbol{x}_B)=0$. Here we build a biquadratic Lagrange polynomial over a $3\times 3$ stencil box composed only by either other ghost points or internal points near $\boldsymbol{x}_G$ (which is included in the stencil).
Let $\xi=\left(\xi_x, \xi_y\right)$ be the distance between the boundary and the ghost point given by
\begin{equation}\label{xi:theta}
    \xi_x = 1-\frac{\abs{\boldsymbol{x}_{B_x} - \boldsymbol{x}_{G_x}}}{h}\qquad \xi_y = 1-\frac{\abs{\boldsymbol{x}_{B_y} - \boldsymbol{x}_{G_y}}}{h},
\end{equation}
with $\boldsymbol{x}_\star=\left(\boldsymbol{x}_{\star_x}, \boldsymbol{x}_{\star_y}\right)$, such that $0\leq \xi_x < 1$ and $0\leq \xi_y < 1$.\\

A Dirichlet-type boundary condition $u(\boldsymbol{x}_B) = g_D$ can be approximated with
\begin{equation}\label{bc:Dirichlet}
    u(\boldsymbol{x}_B) \approx \sum_{i,j=1}^{3}\alpha_i\alpha_j u_{ij} = g_B,
\end{equation}
where $\alpha_i$ are the $2^{\text{nd}}$ order Lagrange interpolation weights given by
\begin{equation}\label{dir:weights}
    \left(\alpha_1, \alpha_2, \alpha_3\right) = \left(\frac{\xi(\xi+1)}{2}, 1-\xi^2, \frac{\xi(\xi-1)}{2}\right).
\end{equation}
$\alpha_i$ represent the weights for the $x$-direction, while $\alpha_j$ are the weights for the $y$-direction, thus they are respectively evaluated with $\xi_x$ and $\xi_y$.\\
On the other hand, a Neumann-type boundary condition $\nabla_n u_{|_{\boldsymbol{x}_B}}=g_N$ is interpolated by
\begin{equation}
    \nabla_n u_{|_{\boldsymbol{x}_B}} \approx n_x\sum_{i,j=1}^{3}\omega_i\alpha_j u_{ij} + n_y\sum_{i,j=1}^3 \alpha_i\omega_j u_{ij} = g_N,
\end{equation}
where $\omega_i$ are the $2^{\text{nd}}$ order Lagrange interpolation weights given by 
\begin{equation}\label{neu:weights}
    \left(\omega_1, \omega_2, \omega_3 \right) = \frac{1}{h}\left(-\frac{1}{2}-\xi,\ 2\xi,\ \frac{1}{2}-\xi\right)
\end{equation}
evaluated in the same way of the weights $\alpha$, and $n_x$, $n_y$ are the normal components.	To simplify, we write
\begin{equation}\label{gh:dirfin}
    \sum_{p=1}^{9} \alpha_{D,p} u_p := \sum_{i,j=1}^{3}\alpha_i\alpha_j u_{ij},
\end{equation}
and
\begin{equation}\label{gh:neufin}
    \sum_{p=1}^{9} \omega_{N,p} u_p := n_x\sum_{i,j=1}^{3}\omega_i\alpha_j u_{ij} + n_y\sum_{i,j=1}^3 \alpha_i\omega_j u_{ij}.
\end{equation}
Then, plugging it into the boundary condition for the humidity equation in the first line of \eqref{model:bc} we get:	
\begin{equation}\label{bound1}
	\left(\frac{n_*}{\tilde{n}}\right)^2 \sum_{p=1}^{9} \omega_{N,p} B_p = K_w \left(\mathcal{E} - \sum_{p=1}^{9} \alpha_{D,p} \theta_p\right),
\end{equation}
where $B_p = B(\theta_p/n_p)$.\\

In case a $3\times3$ stencil of internal and ghost points cannot be constructed (e.g., close to pointly features of the domain), a 9 point stencil and coefficients $\alpha$ and $\omega$ for (\ref{gh:dirfin}) and (\ref{gh:neufin}) can usually be constructed (see \cite{CocoFD}).

We remark that the values of $\theta$ and $c_a$ evaluated at ghost points represent virtual values that satisfy the boundary conditions. They do not represent the real values at the boundary, but they can be recovered easily through interpolation. On the other hand, $n_*$ in \eqref{bound1} represents the porosity at the boundary, which here is always used as the value of the first internal point.

\section{Numerical tests}\label{sec:tests}
In this Section, we aim at simulating the effect of carbonic acid penetrating the porous material using model \eqref{modeleq}. Model parameters are assumed as in Tables \ref{table:stone}, \ref{table:fitting} and \ref{table:kappa}. 

Here, similarly to the approach in \cite{Bracciale}, in order to simulate an accelerated weathering condition we consider a very high value of external humidity $\mathcal{E} \simeq 0.00138$ that we calibrated numerically. 
More generally, in order to make simulations with environmental data we can obtain the value of external humidity by the formula in \cite{HyperPhysics}:
\begin{equation}
	u_R \left( \left(5.018 + 0.32321 \cdot T + 8.1847\cdot 10^{-3} T^2 + 3.1243\cdot 10^{-4} T^3\right) \cdot 10^{-6}\right),
\end{equation}
where $T$ and $u_R$, are, respectively, the ambient temperature and the  percentage relative humidity of the environment, that can be time varying or averaged over time. 

For the external value of $CO_2$, we assume the mean value $\mathcal{C}=$ 5.9e-07 $g/cm^3$ reported in Figure 8 of the paper \cite{bretticeseri}.

We performed the following three experiments:
\paragraph{1D-line.} The computational domain is $\left[0, 5.5\right]$, discretized on a one-dimensional grid with 100 points, which contains the rock sample as the segment $\left[0.25, 5.25\right]$ cm. For the time discretization with $t_{\tau}:=\tau \Delta t$  we choose $\Delta t = 1$ second.

\paragraph{2D-rectangle.} The computational domain is the square $\left[0, 5.5\right]^2$, with the rock sample identified as the rectangle $\left[0.25, 5.25\right] \times \left[0.75, 4.75\right]$ cm. Then the domain is discretized with $100^{2}$ points on a uniform two-dimensional grid, with $\Delta t=0.1$ seconds.

\paragraph{2D-brick.} The computational domain and the brick sample are the same as the previous experiment, however we removed two disks from the bottom edge
\begin{displaymath}
	C_1 : (x-1.5)^2 + (y-1)^2 < 0.5^2 \qquad C_2 : (x-4)^2 + (y-1)^2 < 0.5^2,
\end{displaymath}
and added a small bump at the top edge using the disk
\begin{displaymath}
	C_3 : (x-2)^2 + (y-4.75)^2 \leq 0.25^2.
\end{displaymath}
\begin{figure}
    \centering
    \includegraphics[width=0.7\linewidth]{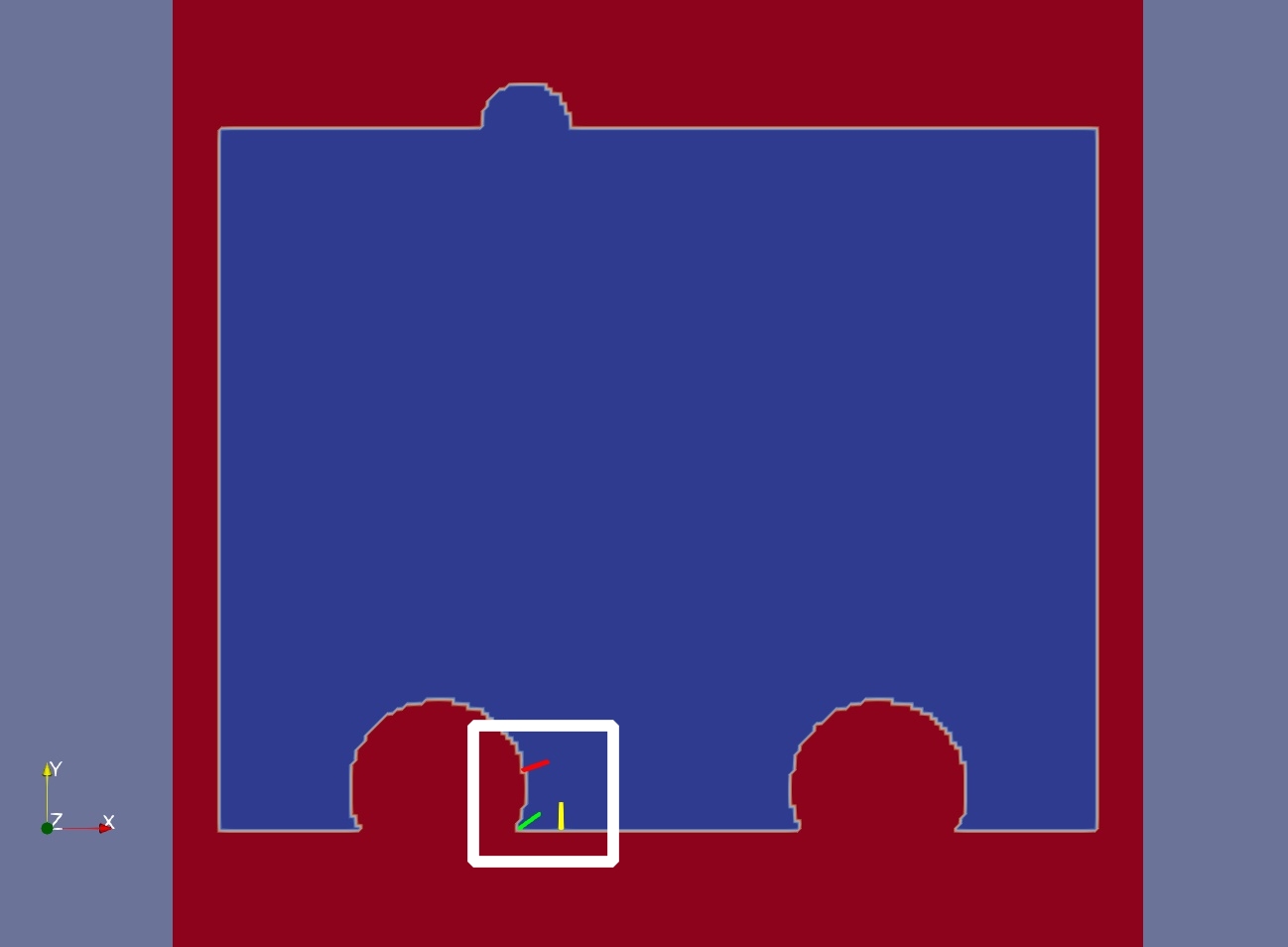}
    \caption{Shape of the domain of the 2D-brick experiment.}
    \label{fig:brickshape}
\end{figure}
In Figure \ref{fig:brickshape} we show the shape of this domain, highlighting with a white box the area where we are going to display the simulation results.\\
In order to increase the resolution of the curved parts of the boundary, we discretized the domain with $200^2$ points on a uniform two-dimensional grid. For increased robustness, we applied a variable timestepping mechanism. We start the simulation with $\Delta t$ equal to half of the spatial step and when the nonlinear solver converges in less than 6 iterations, we try doubling the next timestep, up to a maximum of $\Delta t = 6$ seconds.\\

\begin{figure}
	\centering
	\includegraphics[width=0.44\linewidth]{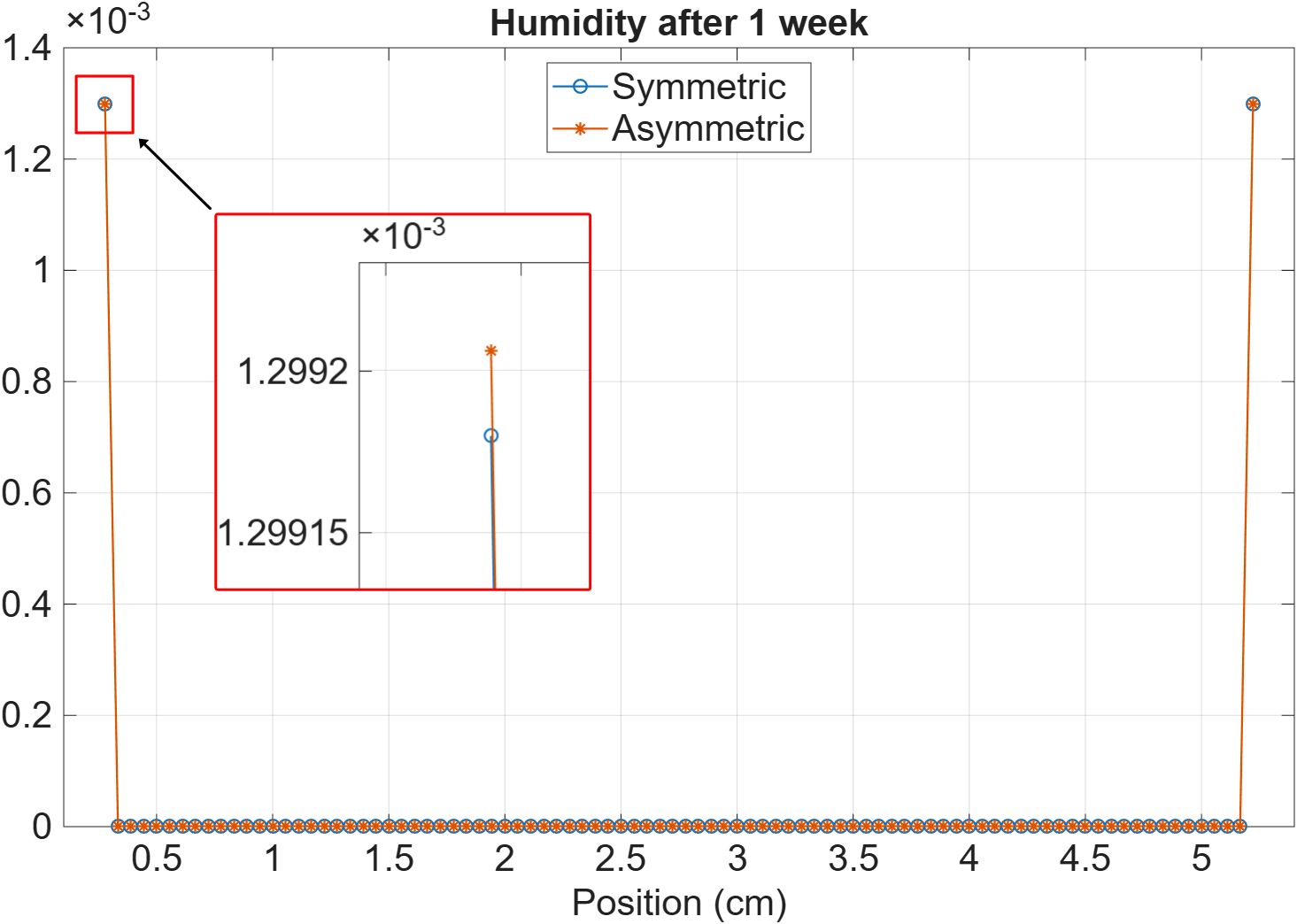}
	\includegraphics[width=0.44\linewidth]{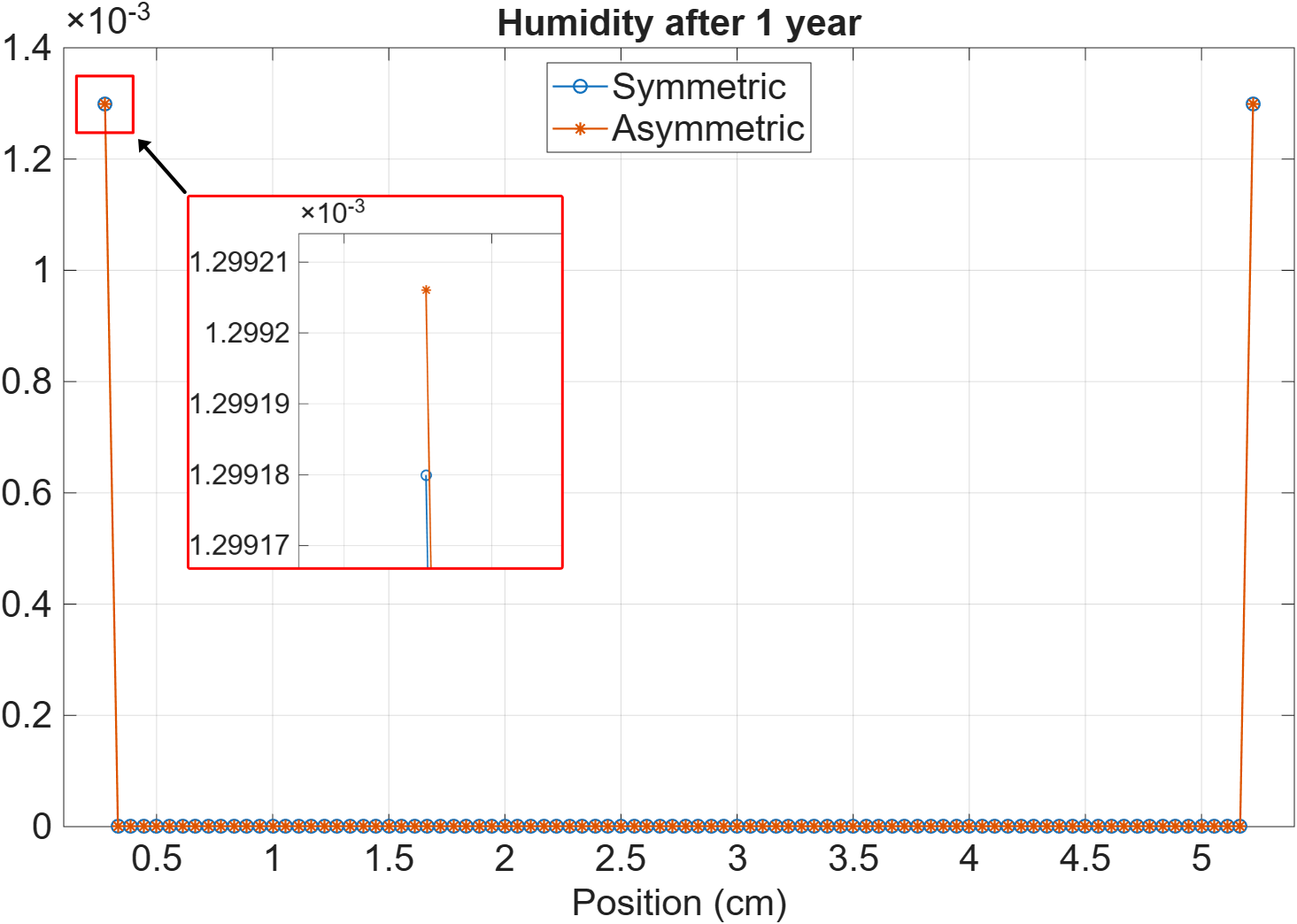}
	\caption{Plots of the humidity of 1D simulation at two different timestamps. Left panel: output of symmetric and non-symmetric models after 1 week. Right panel: output of symmetric and non-symmetric models after  1 year.}
	\label{fig:humiditycomp}
\end{figure}

First we depict in Figure \ref{fig:humiditycomp} the plots obtained in 1D for the humidity $\theta$ with both symmetric (blue line-points) and asymmetric model (red line-points), to show that they produce results that are almost identical. Since the asymmetric absorption function $B'_{kP}$ has parameters that make it sufficiently close to the symmetric one in terms of diffusivity coefficient and shape, it is not surprising that the numerical solutions do not differ greatly from each other. The same can be said for the positions of the boundary, as it will be shown in Fig. \ref{fig:boundary}.
Then, in the next pictures we report the plots obtained with asymmetric model, since it is more flexible due to its dependence on a greater number of parameters.\\

The parameter set $P:=\left(K_{w},K_{a},K_{n},n_{max}, \mathcal{E}\right)\in\mathbb{R^+}$ reported in Table \ref{table:kappa} is calibrated in 1D model by solving the minimization problem  
\begin{eqnarray}
	\min_{P} F\left(P\right):=\frac{1}{5}\displaystyle \sum_{i=1}^{5}\left(\frac{\sigma^{num}\left(P,t_{i}\right)-\sigma(t_i)}{\sigma(t_i)}\right)^{2},
\end{eqnarray}

at  $t_i=1,\dots,5$, with $t_i =7, 30, 90, 180, 365$ days,
with cost function $F(P)$, which is the $L^{2}$-norm error between the numerical values of recession front $\sigma^{num}(P,t_i)$ and the experimental ones defined as the linear function $\sigma_{exp}(t_i)= r t_i$ with $r$ the recession rate reported in Table \ref{table:stone}. Note that the numerical values are computed by the forward model (\ref{modeleq}) for a given parameter set $P$ as:
$\sigma^{num}\left(P,t_{i}\right):=\min_{x} \left\lbrace |x-x_0| s.t. \ \varphi(x,t_{i})=0\right\rbrace$ with $x_0$ the initial left boundary of the object and $\varphi(x,t)$ defined through the level set (\ref{levelset}).
\paragraph{Dataset.} The recession rate $r=2 \ mm/ 100 \ years$ is not used as a single scalar data point to identify the five model parameters of $P$. Since experimental evidence indicates that, once steady-state erosion due to weathering is established, the stone surface recedes at an essentially constant velocity \cite{brimblecombe,Hoke2002,Hoke2004}, the recession front position is known for every time through the linear law $\sigma_{exp}(t_{i}) = r t_{i}$, which we sample at sufficient time points $t_{i}$. Physically, this selects parameter sets producing a  regime, in which the penetration profile of $c_a$ has reached its travelling-wave shape and the free boundary $\Gamma(t)=\{\varphi(\cdot,t)=0\}$ advances
at constant speed, and rejects those that reproduce the correct material loss after one year through a spurious transient, such as an incubation stage followed by an accelerating front. A genuine enlargement of the information content requires
observables of a different nature, e.g.\ the incubation time or the thickness of the
dissolution layer, which we discuss in Section~\ref{sec:conclusions} as future work.
\\
\\
Using simulated annealing as a global nonlinear optimization method, we solve the inverse problem to calibrate the optimal parameter set $P^{opt}$ reported in Table \ref{table:kappa} with recession rate $r^{num}=2.0665 \ mm/100 \ years$ compared to the recession rate  $r=2 \ mm/100 \ years$ from literature (see Table \ref{table:stone}) with cost function value $F(P^{opt})= 2.149 \cdot 10^{-2}$. 
The calibrated parameters are used in the following numerical tests.\\

\begin{figure}
	\centering
	\includegraphics[width=0.48\linewidth]{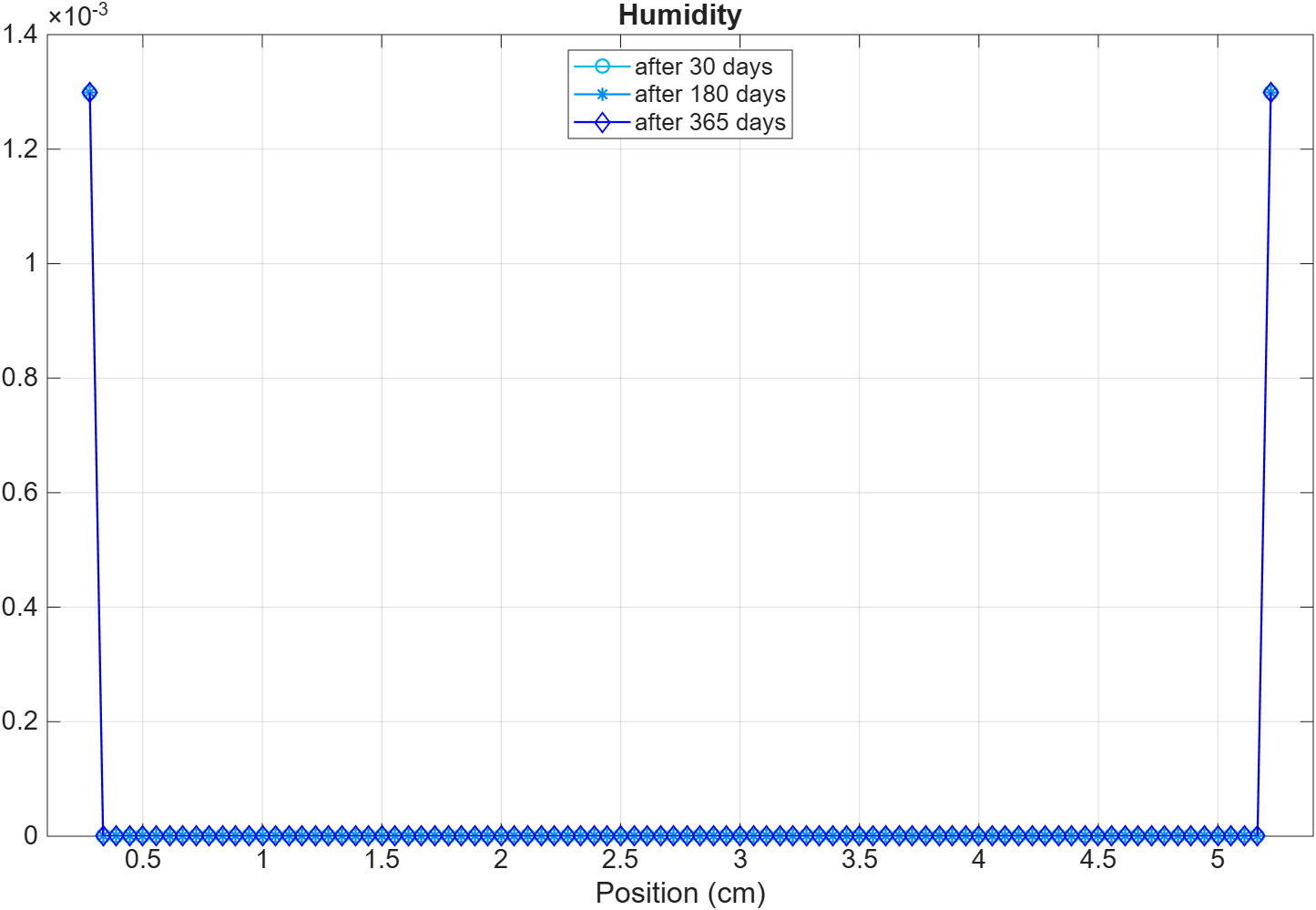}
	\includegraphics[width=0.48\linewidth]{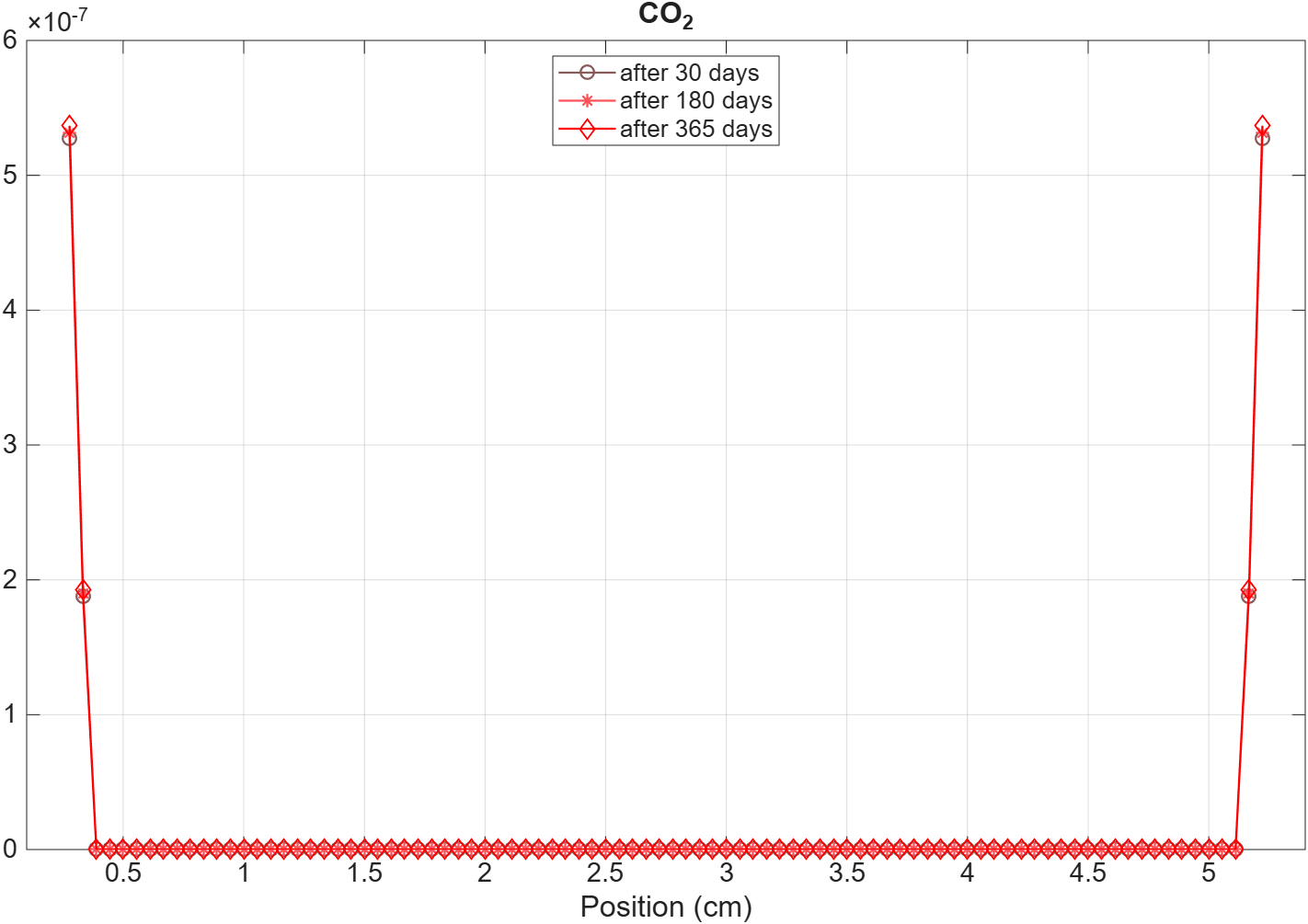}\\
	\vspace{0.2cm}
	\includegraphics[width=0.48\linewidth]{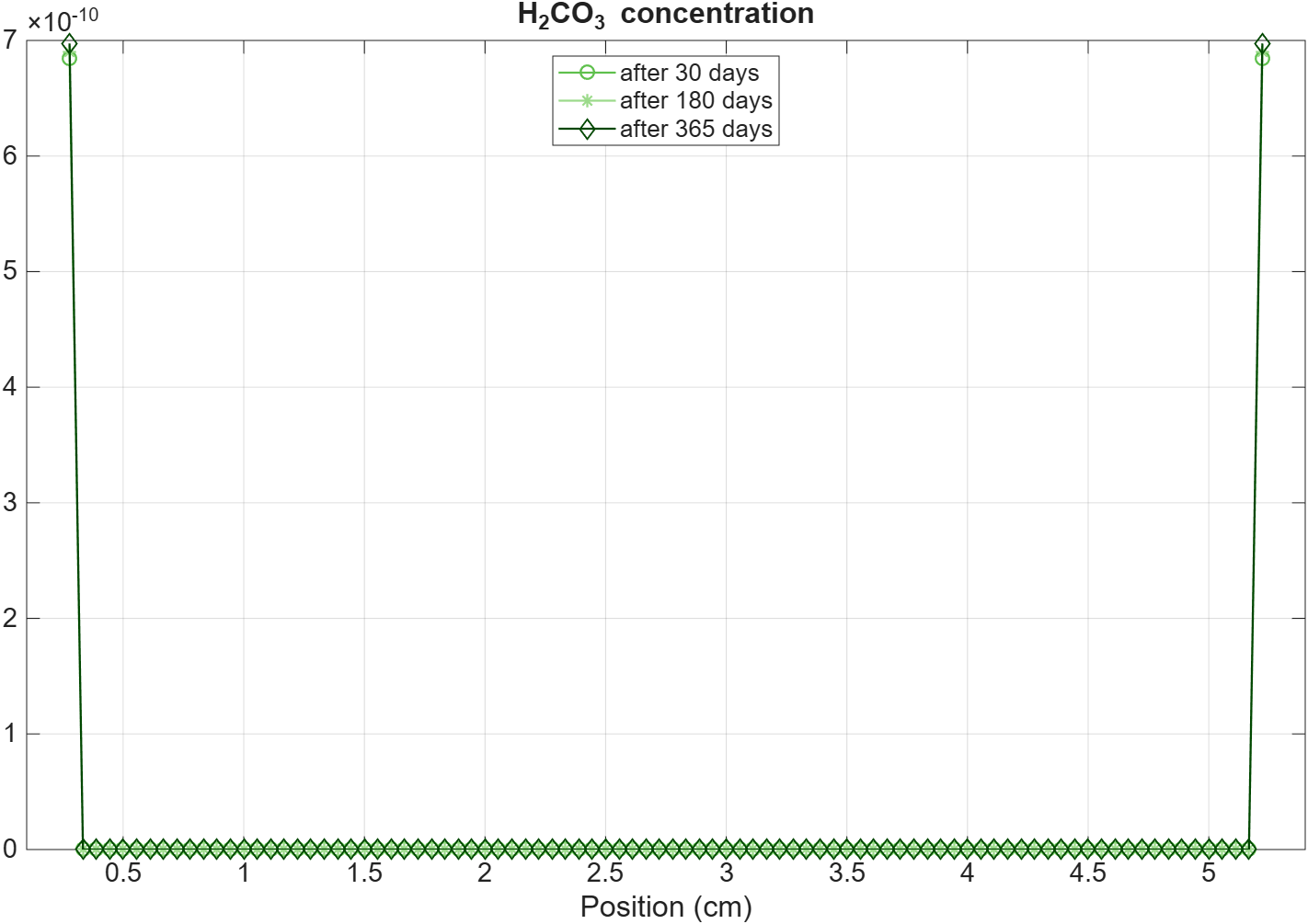}
	\includegraphics[width=0.48\linewidth]{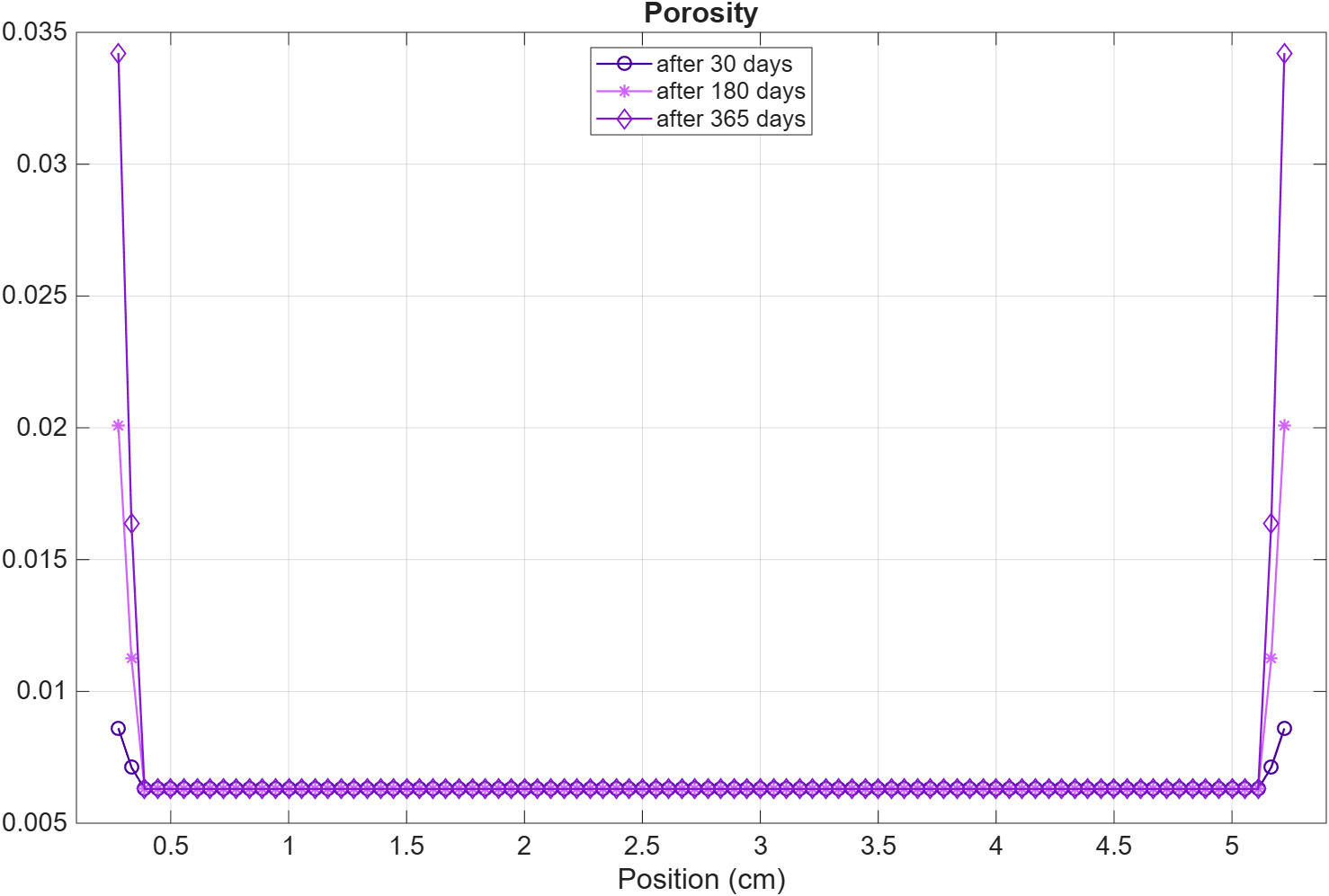}
	\caption{Plots of the variables of 1D simulation at three different timestamps: one week, 180 days and one year. The simulation is performed with the asymmetric $B'_{kP}(s)$. From top left: humidity, $CO_2$, carbonic acid and porosity.}
	\label{fig:asymSim}
\end{figure}

\begin{figure}
	\centering
	\includegraphics[width=0.48\linewidth]{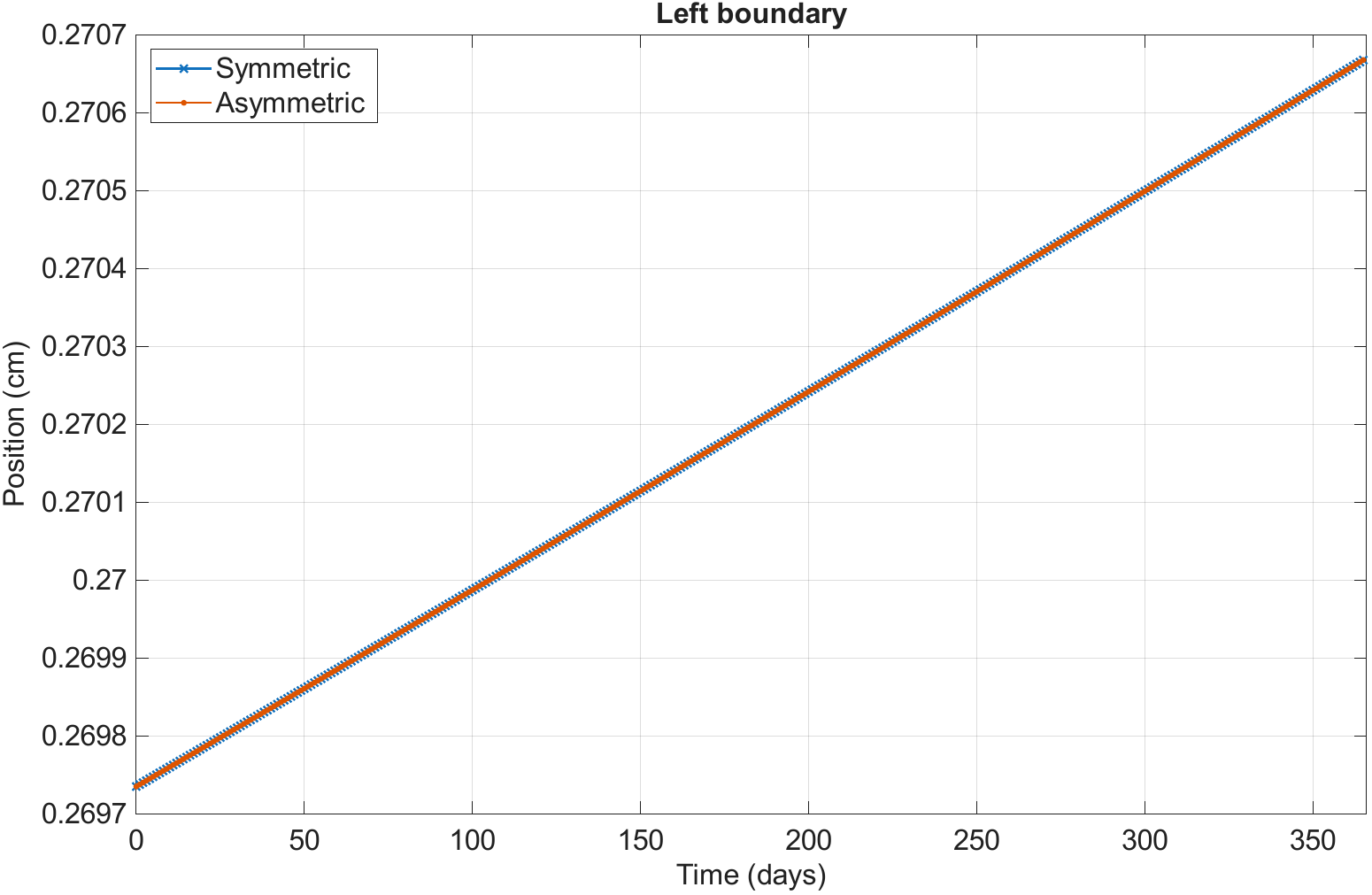}
	\includegraphics[width=0.48\linewidth]{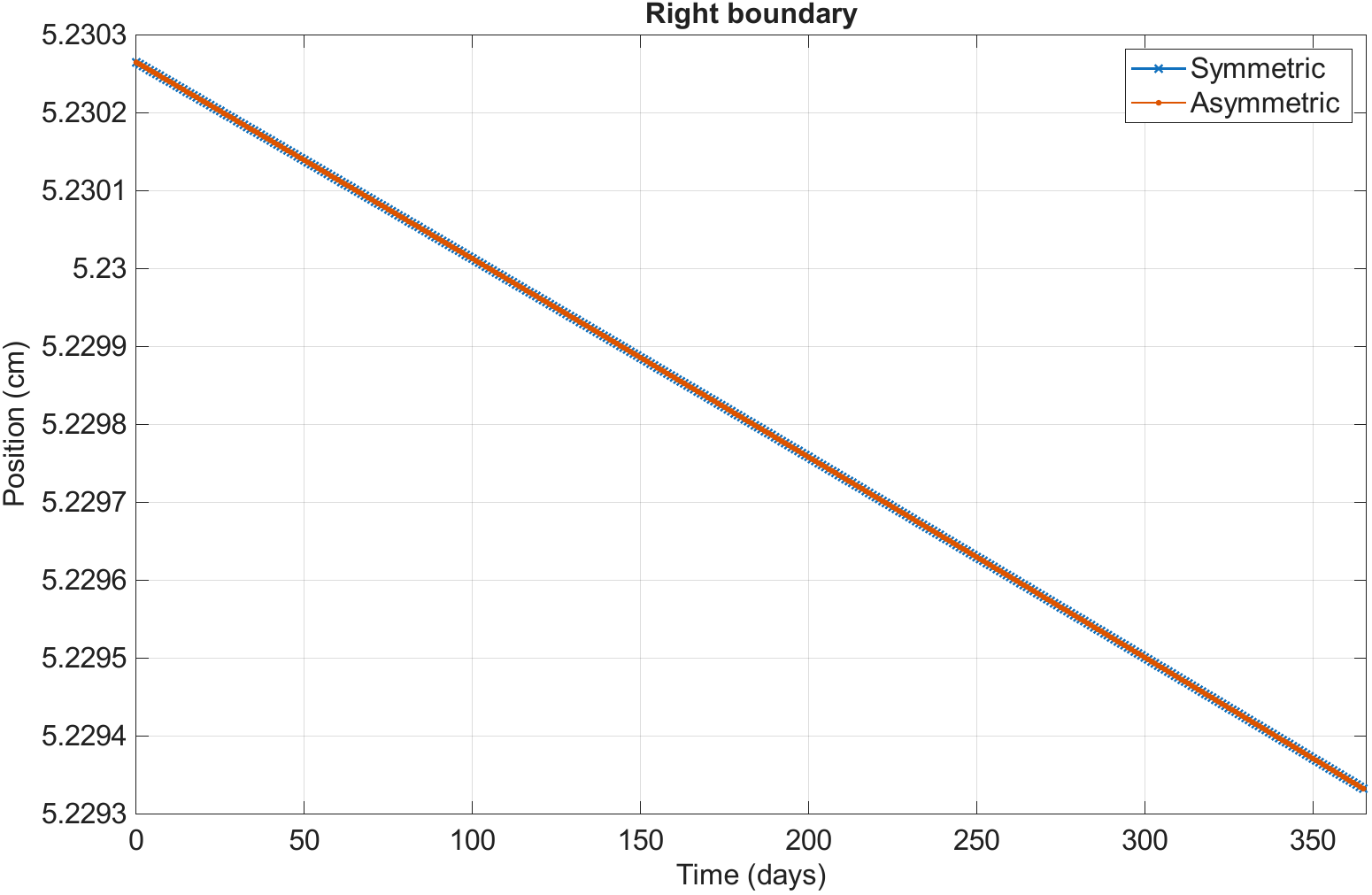}
	\caption{Plots of the moving boundaries. Left panel: left side of the specimen. Right panel: right side of the specimen.}
	\label{fig:boundary}
\end{figure}

In Figure \ref{fig:asymSim} we show the numerical results obtained with a simulation of 1 year with non-symmetrical $B_{kP}'(s)$ and parameters as in Table \ref{table:fitting}. As expected the numerical scheme reaches the steady state solutions for both humidity (top left) and dissoluted carbonic acid (bottom left) quite fast, while the porosity (bottom right) is increasing until material loss occurs. Finally, on the top right panel of Fig. \ref{fig:asymSim} we show the $CO_2$ level penetrating into the porous matrix across time.

Figure \ref{fig:boundary} shows the positions of the edges over time; here too, the results are shown for both symmetrical and asymmetrical $B(s)$. However, as it can be seen from the results in Table \ref{table:erosion}, the total boundary movement is different between the two simulations. This means that if the model parameters in the asymmetrical absorption function were calibrated against a richer dataset, thus making the shape of $B'_{kP}$ less adherent to that of the symmetrical one, we would have observed a greater difference.

\begin{table}
	\centering
	\begin{tabular}{|c|c|c|}
		\toprule
		$B(s)$ & Erosion of the left edge [cm] & Erosion of the right edge [cm]\\
		\midrule
		Symmetric & 0.0206677108748 & 0.0206677108748 \\
		Asymmetric & 0.0206677027187 & 0.0206677027187 \\
		\midrule
		Difference & $8.15 \cdot 10^{-9}$ & $8.15 \cdot 10^{-9}$\\ 
		\midrule
		Recession rate & 6.5528$\cdot 10^{-10} cm \ s^{-1}$& 6.5528$\cdot 10^{-10}cm \ s^{-1}$ \\
		\bottomrule
	\end{tabular}
	\caption{Eroded space after 1 year of simulation.}
	\label{table:erosion}
\end{table}

\begin{figure}
	\centering
	\includegraphics[width=0.7\linewidth]{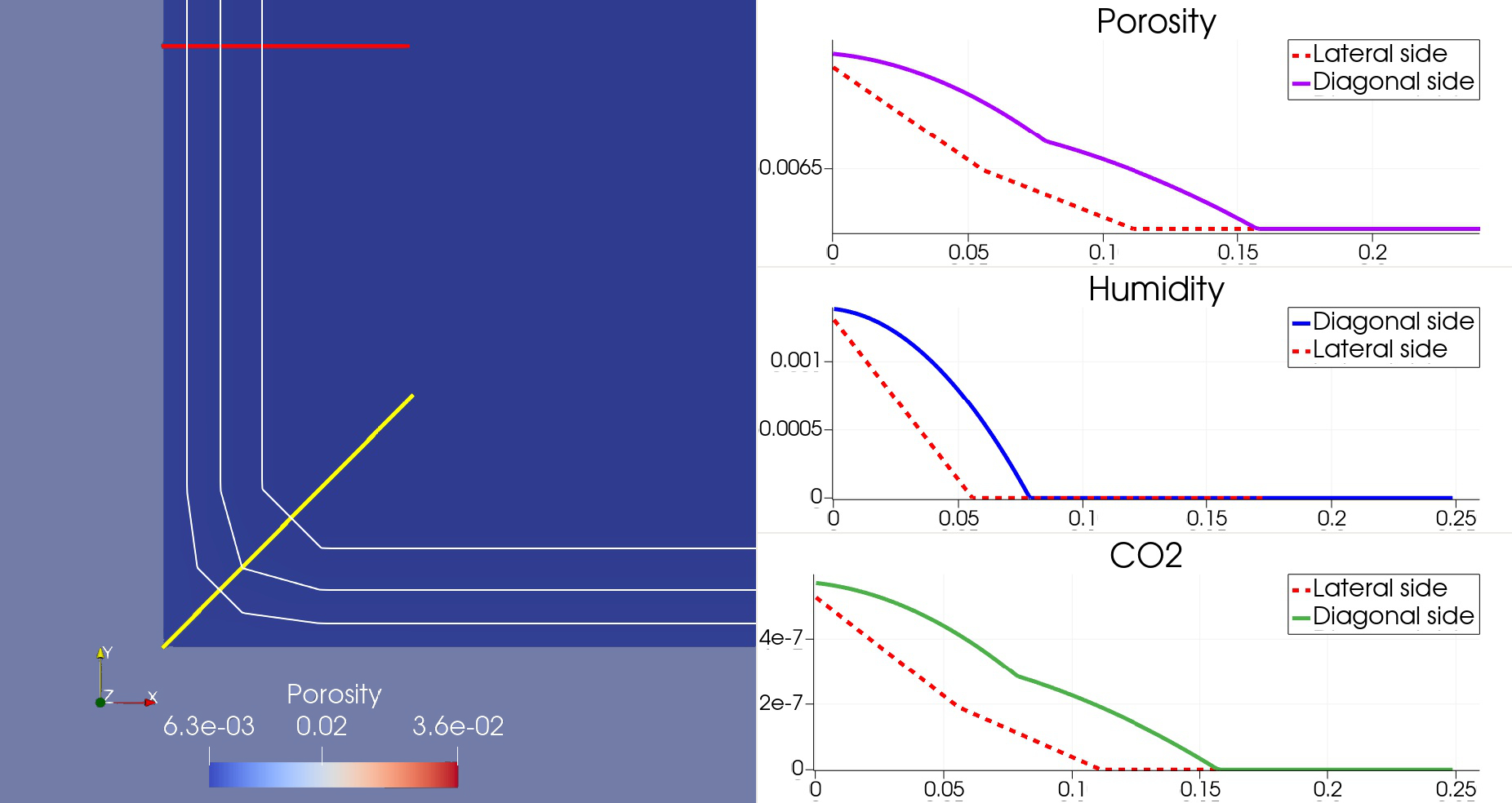}\\
	\vspace{1mm}
	\includegraphics[width=0.7\linewidth]{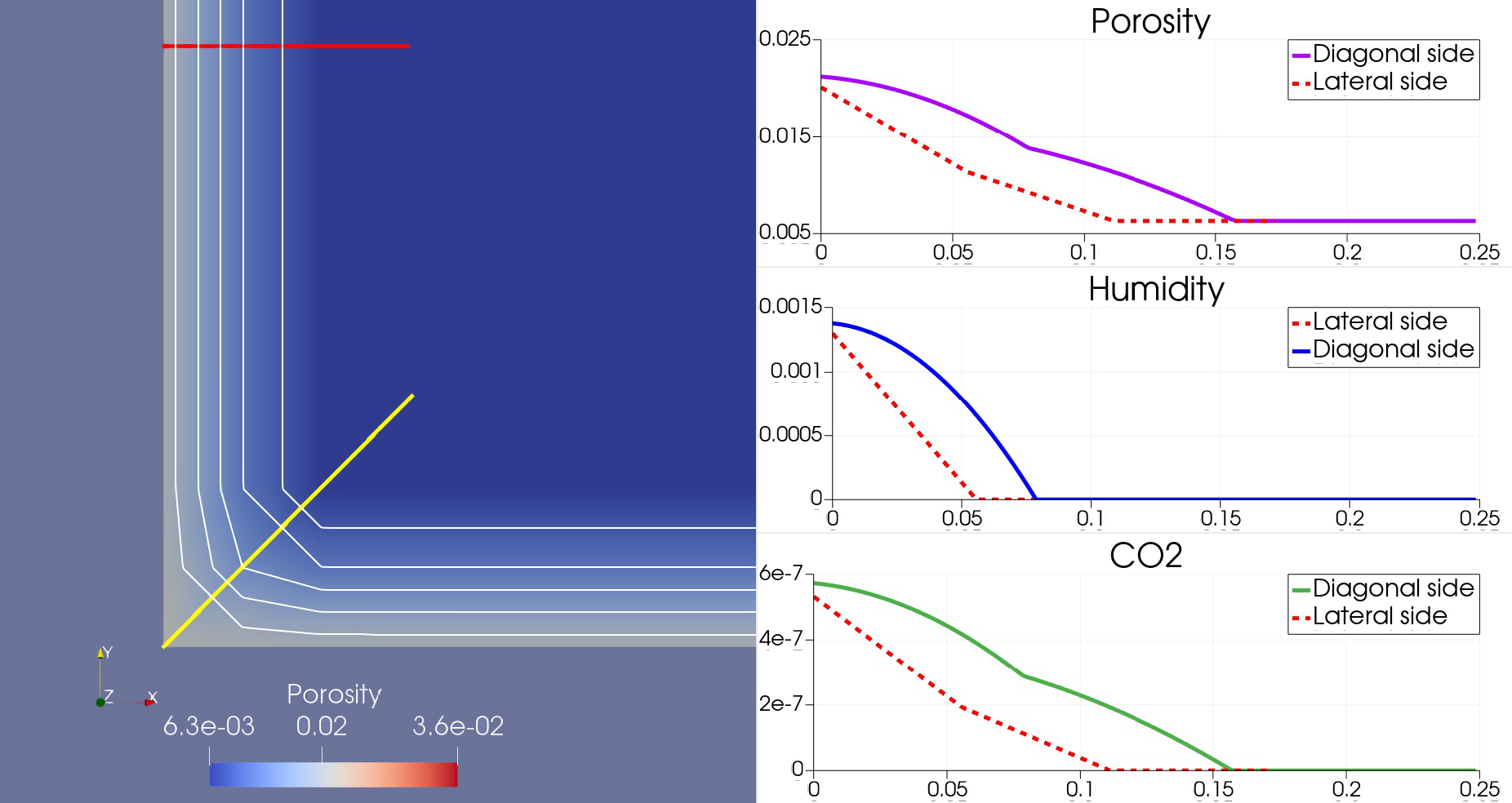}\\
	\vspace{1mm}
	\includegraphics[width=0.7\linewidth]{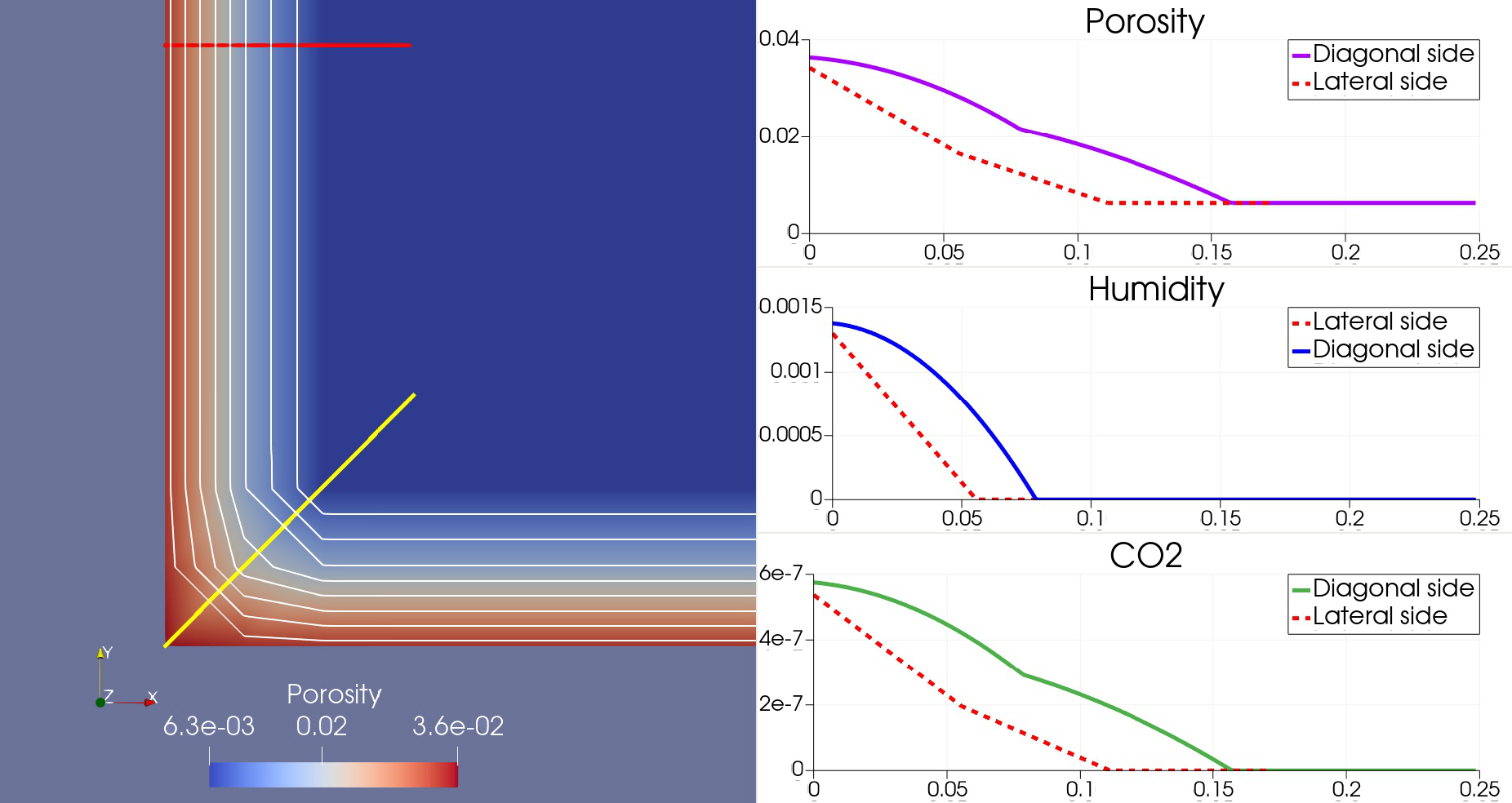}
	\caption{Plots of the 2D-rectangle simulation over one year. From above: simulation after one week, after 180 days, after 1 year. On the right side there are the plots of the variables of the system along two different direction. The solid line is along the diagonal side of the domain, while the red-dashed is along the lateral side.}
	\label{fig:2Dsim}
\end{figure}

\begin{figure}
	\centering
	\includegraphics[width=0.8\linewidth]{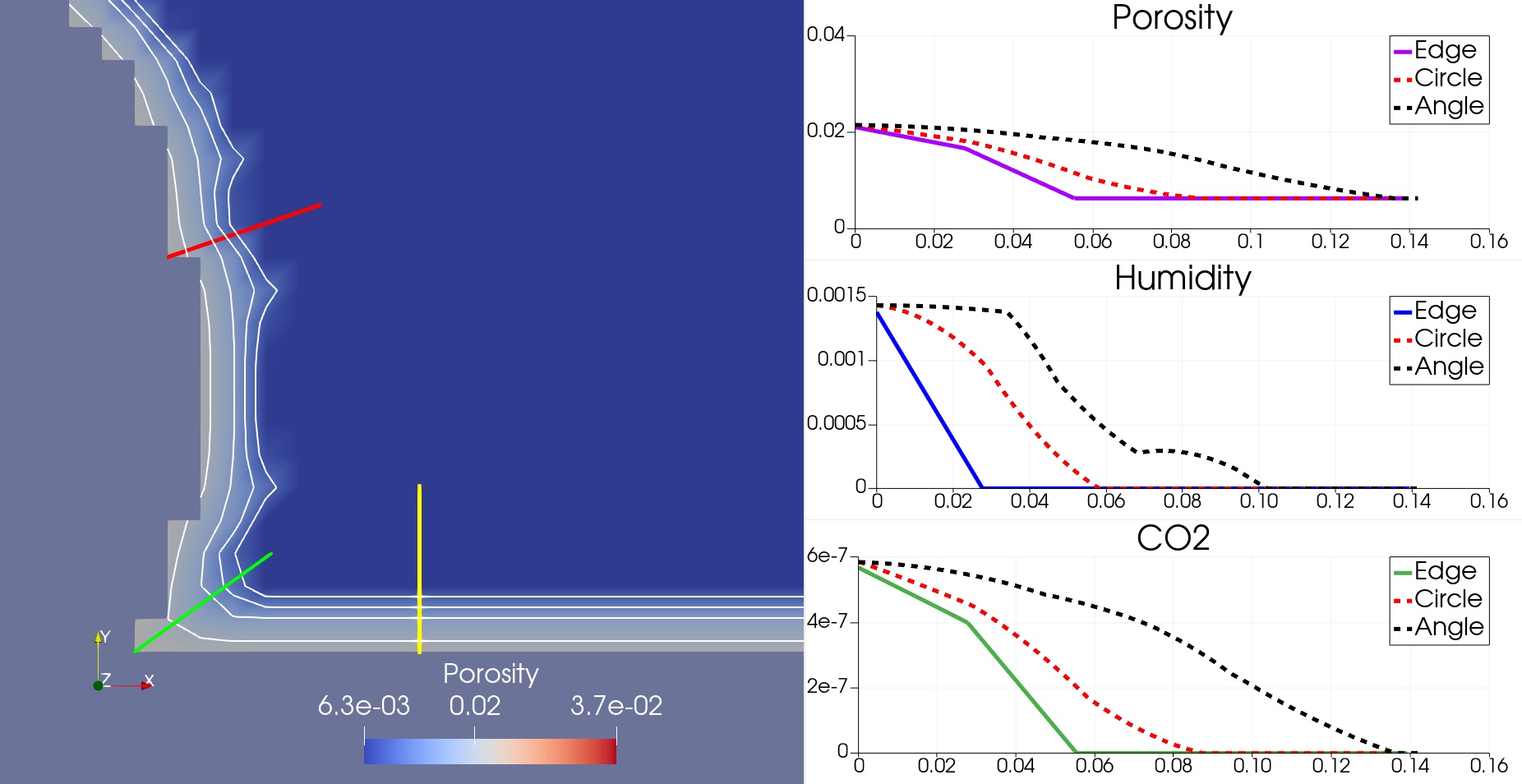}\\
	\vspace{1mm}
	\includegraphics[width=0.8\linewidth]{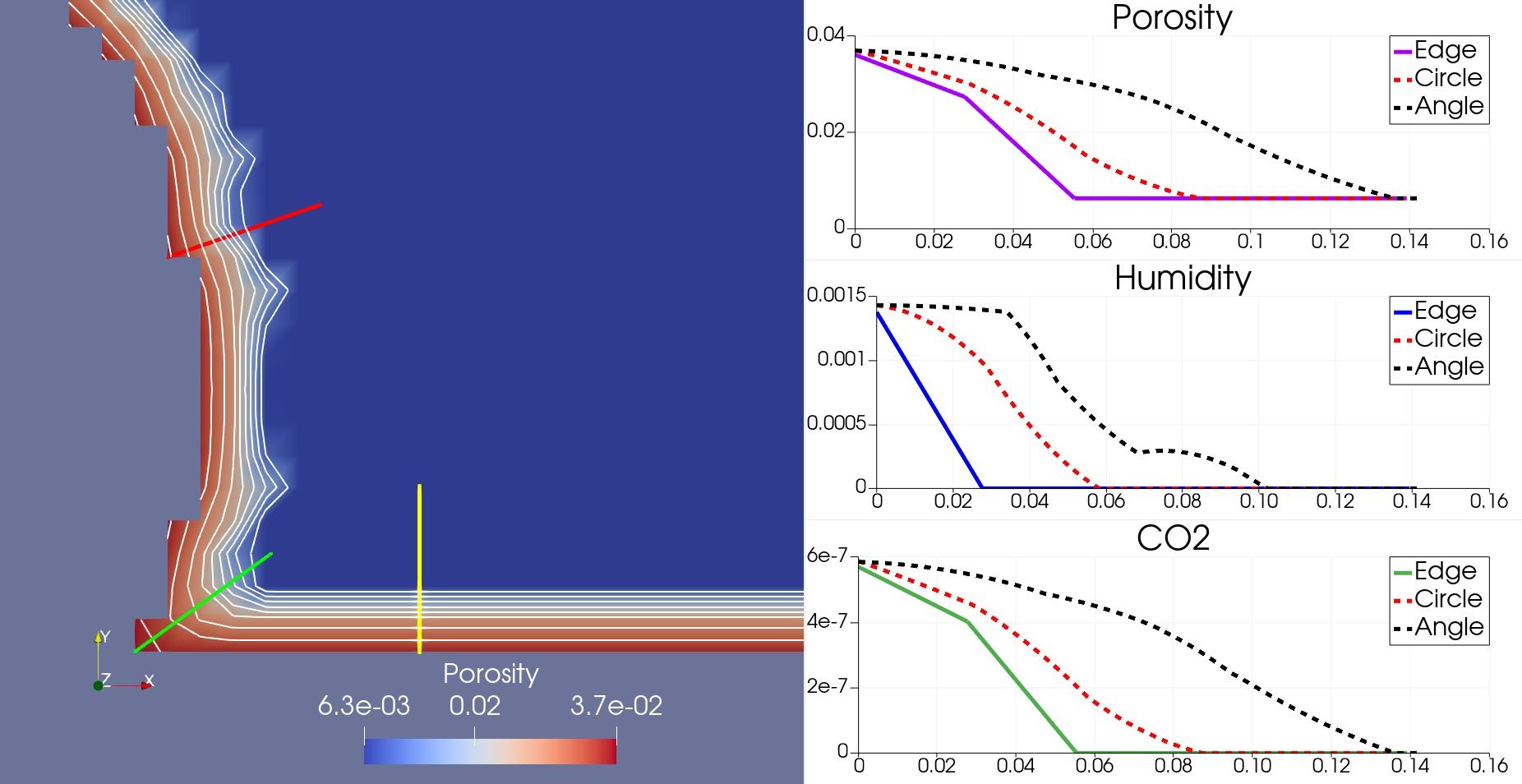}
	\caption{Plots of the 2D-brick experiment over one year. From above: simulation after 180 days and after 1 year. On the right side there are the plots of the variables of the system along three different directions. The solid line is along the edge of the domain, while the red-dashed is along the normal direction of the circle, and the black-dashed is along the normal direction from the angle.}
	\label{fig:2DLego}
\end{figure}

\begin{figure}
    \centering
    \begin{subfigure}{0.45\linewidth}
        \centering
        \includegraphics[width=\linewidth]{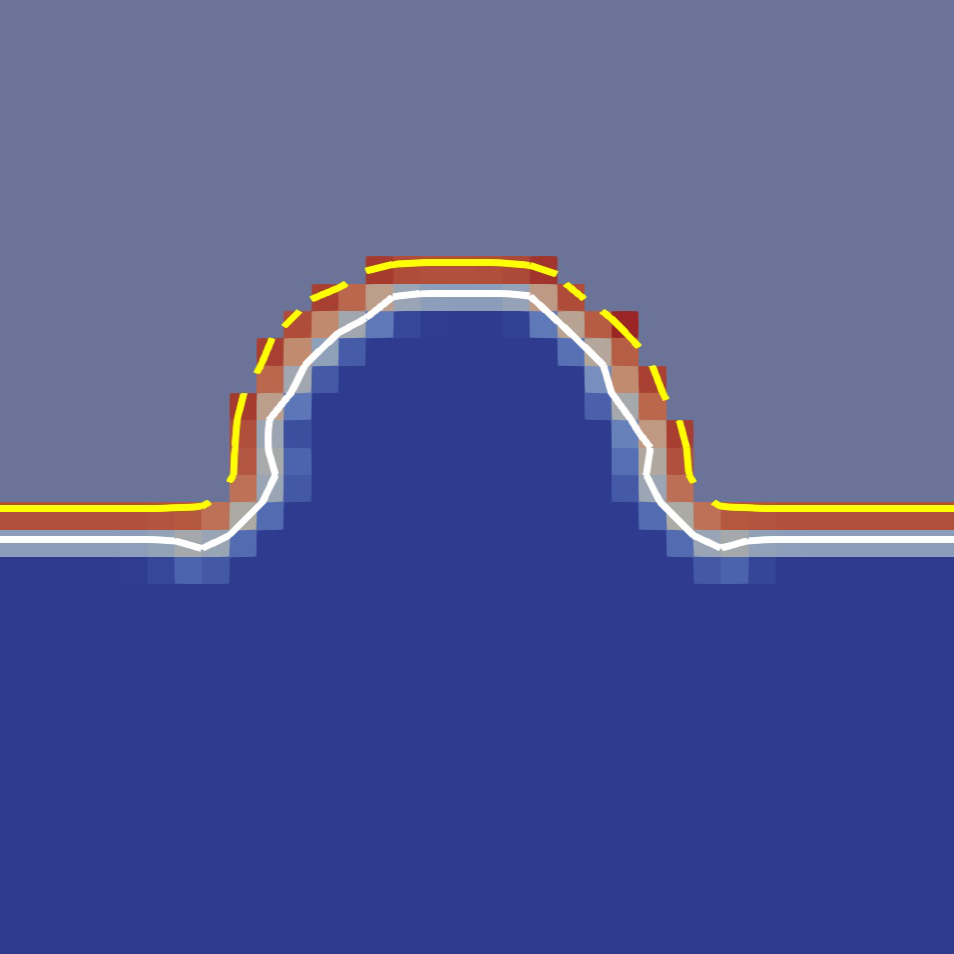}
        \caption{Top edge bump}
    \end{subfigure}
    \hspace{1em}
    \begin{subfigure}{0.45\linewidth}
        \centering
        \includegraphics[width=\linewidth]{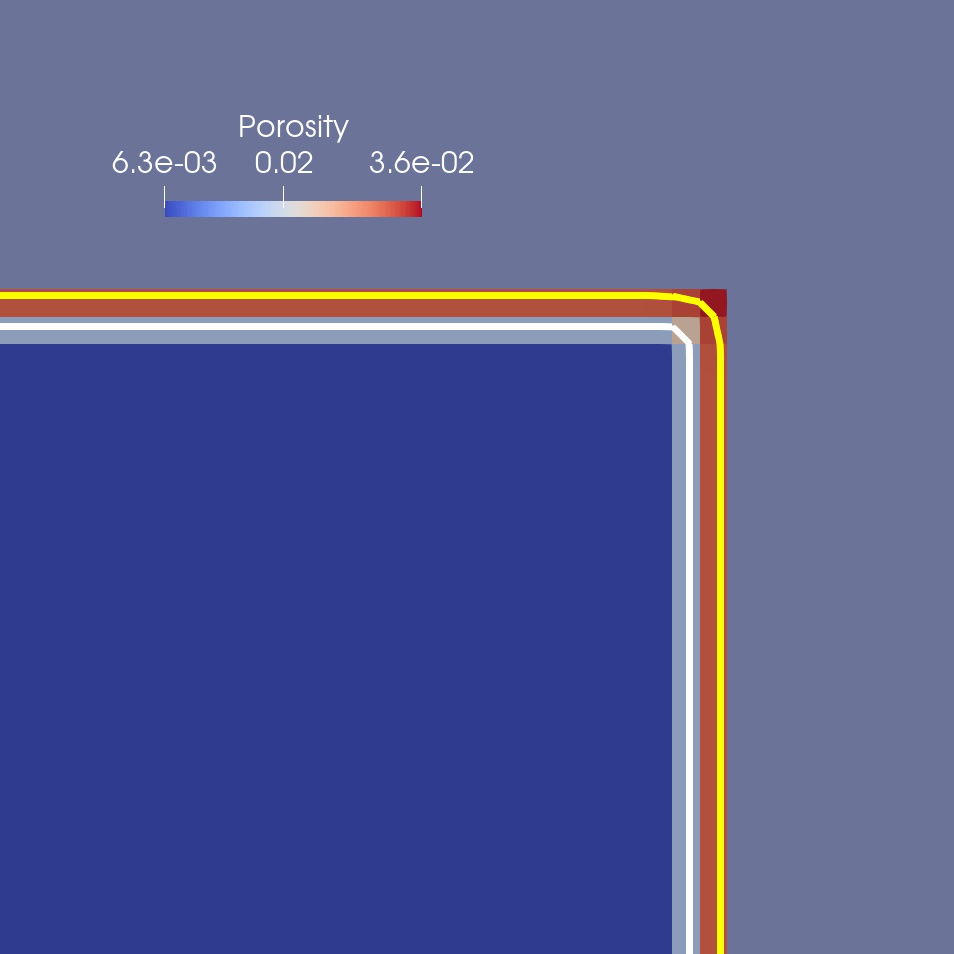}
        \caption{Top right corner}
    \end{subfigure}
    \begin{subfigure}{0.45\linewidth}
        \centering
        \includegraphics[width=\linewidth]{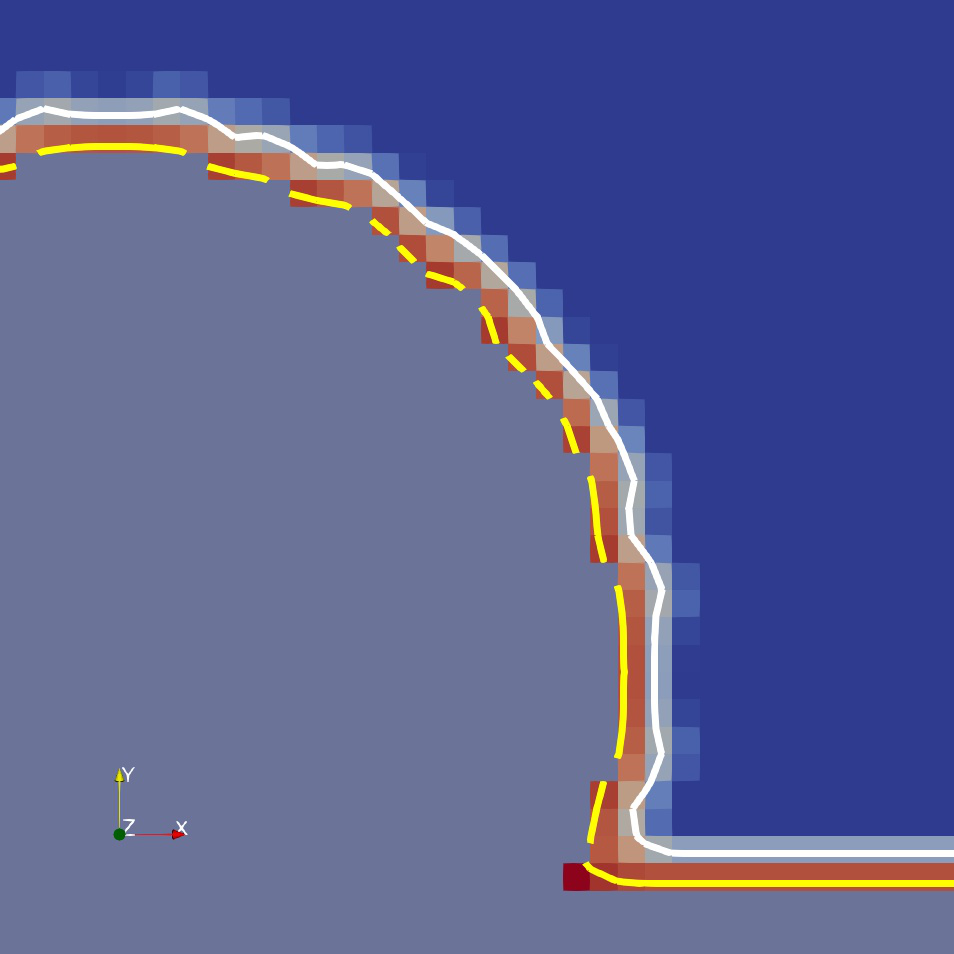}
        \caption{Bottom edge left hole}
    \end{subfigure}
    \hspace{1em}
    \begin{subfigure}{0.45\linewidth}
        \centering
        \includegraphics[width=\linewidth]{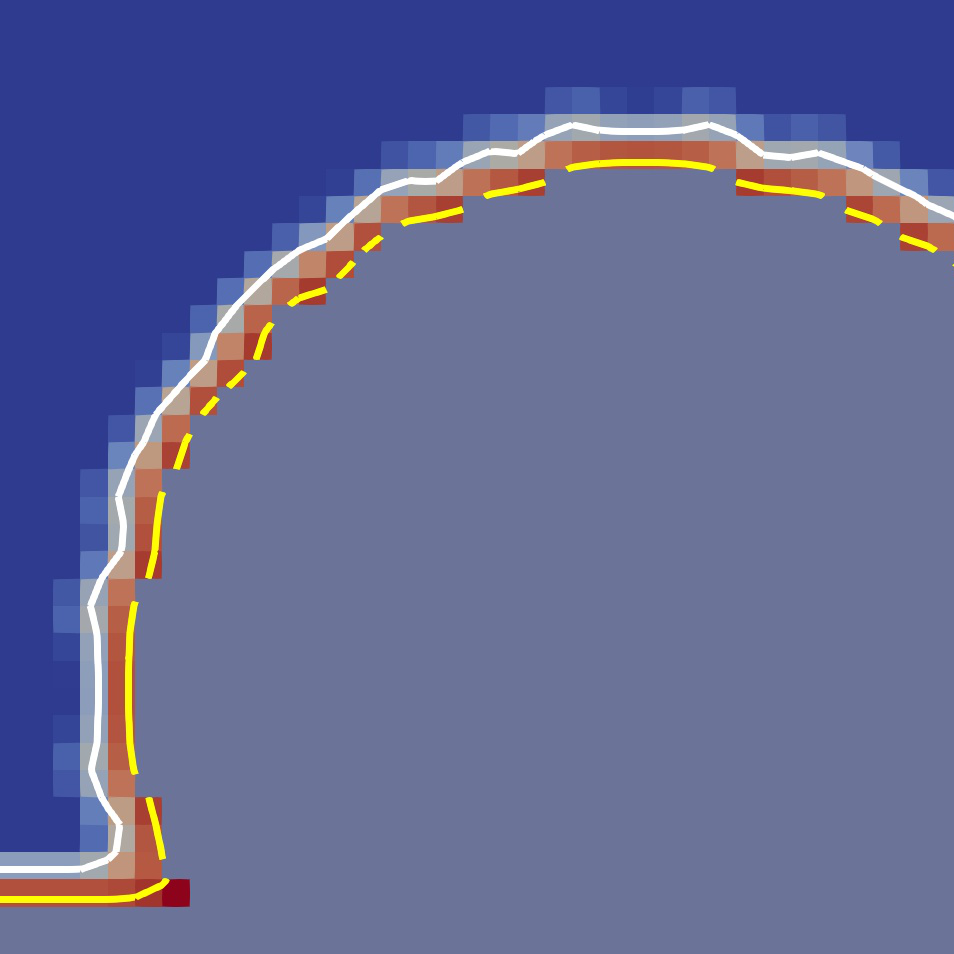}
        \caption{Bottom edge right hole}
    \end{subfigure}
    \caption{Contour of the porosity value $n=0.02$ for the 2D-brick experiment in four different locations. The yellow line represents the position after half year, while in white after 1 year.}
    \label{fig:recessionParagon}
\end{figure}

In Figure \ref{fig:2Dsim} and \ref{fig:2DLego} we show the results of numerical simulation in the 2D settings after 1 year. The parameters are the same as in the one-dimensional case. Figure \ref{fig:2Dsim} shows a detail of the bottom-left corner of the speciment, while Figure \ref{fig:2DLego} shows the detail in the white box from Figure \ref{fig:brickshape}. As before, the porosity shows an increasing behavior. Moreover, in Figure \ref{fig:recessionParagon} we display the contour of the porosity level $n=0.02$ for the 2D-brick experiment at two specific timesteps (180 days in yellow, 1 year in white).

\paragraph{Order of the scheme}
\begin{table}
	\centering
	\begin{tabular}{|c|c|c|cc|cc|}
		\toprule
		N & $\Delta x$ & $\Delta t$ & Err($\theta$) & order & Err($c_a$) & order\\
		\midrule
		100 & $5.50 \cdot 10^{-2}$ & 1 & 0.00792 & -- 	  & $3.57 \cdot 10^{-6}$ & --\\
		200 & $2.75 \cdot 10^{-2}$ & 0.5 & 0.00809 & 1.022 & $4.08 \cdot 10^{-6}$ & 1.143\\
		400 & $1.38 \cdot 10^{-2}$ & 0.25 & 0.00871 & 1.076 & $4.98 \cdot 10^{-6}$ & 1.220\\
		800 & $6.88 \cdot 10^{-3}$ & 0.125 & 0.00984 & 1.130 & $5.49 \cdot 10^{-6}$ & 1.103\\
		1600 & $3.44 \cdot 10^{-3}$ & 0.0625 & 0.00733 & 0.744 & $4.34 \cdot 10^{-6}$ & 0.791\\
		3200 & $1.72 \cdot 10^{-3}$ & 0.03125 & 0.00777 & 1.060 & $4.37  \cdot 10^{-6}$ & 1.008\\
		\bottomrule
	\end{tabular}
	\caption{Discrete 2-norm and order of accuracy of the scheme in 1D.}
	\label{table:order}
\end{table}
In Table \ref{table:order} we record the order of the numerical scheme. As expected, since we are using a first order scheme in time and first order in space, we get a global numerical order of 1.

\paragraph{Catastrophic condition.}
\begin{figure}
	\centering
	\includegraphics[width=0.49\linewidth]{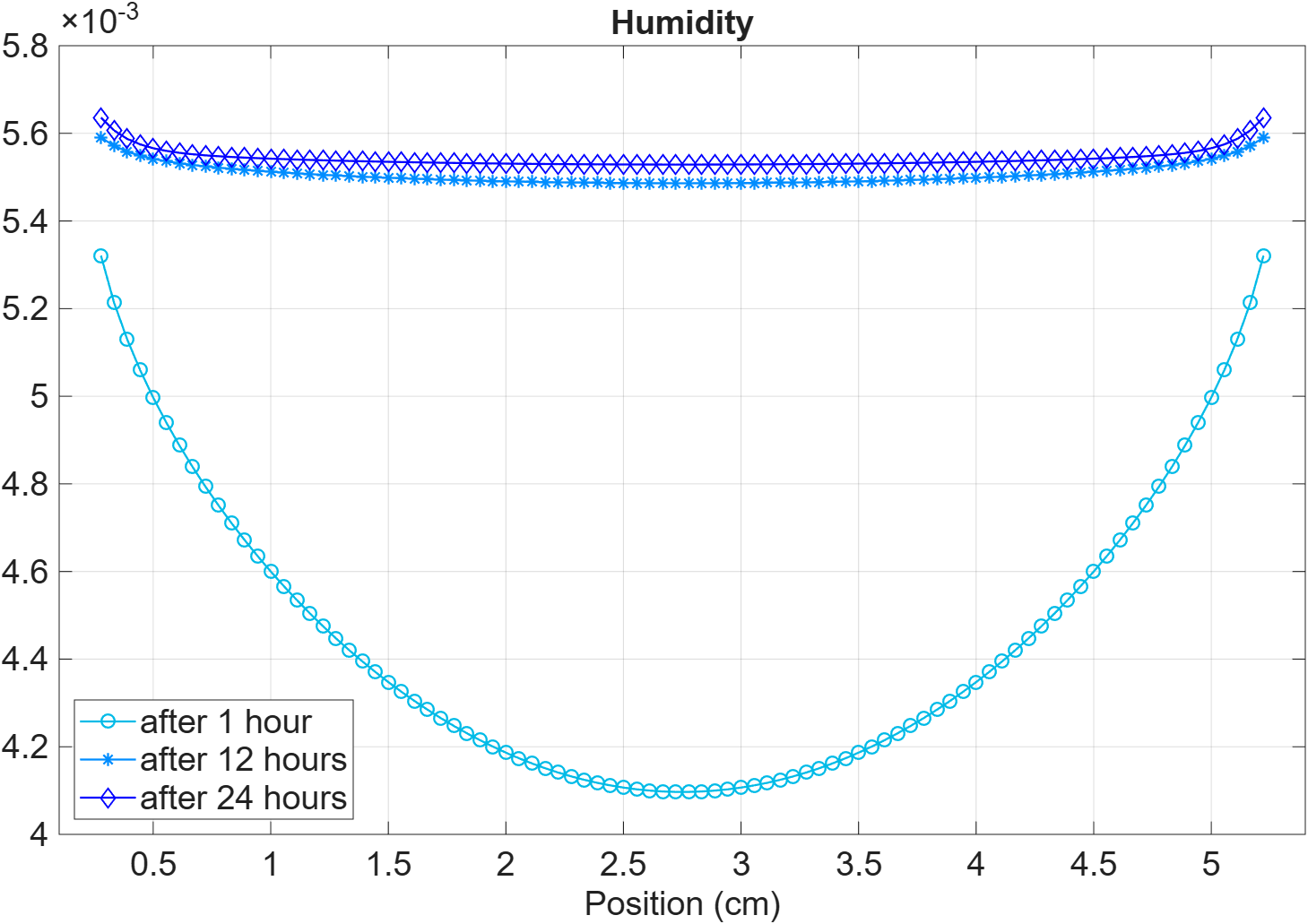}
	\includegraphics[width=0.49\linewidth]{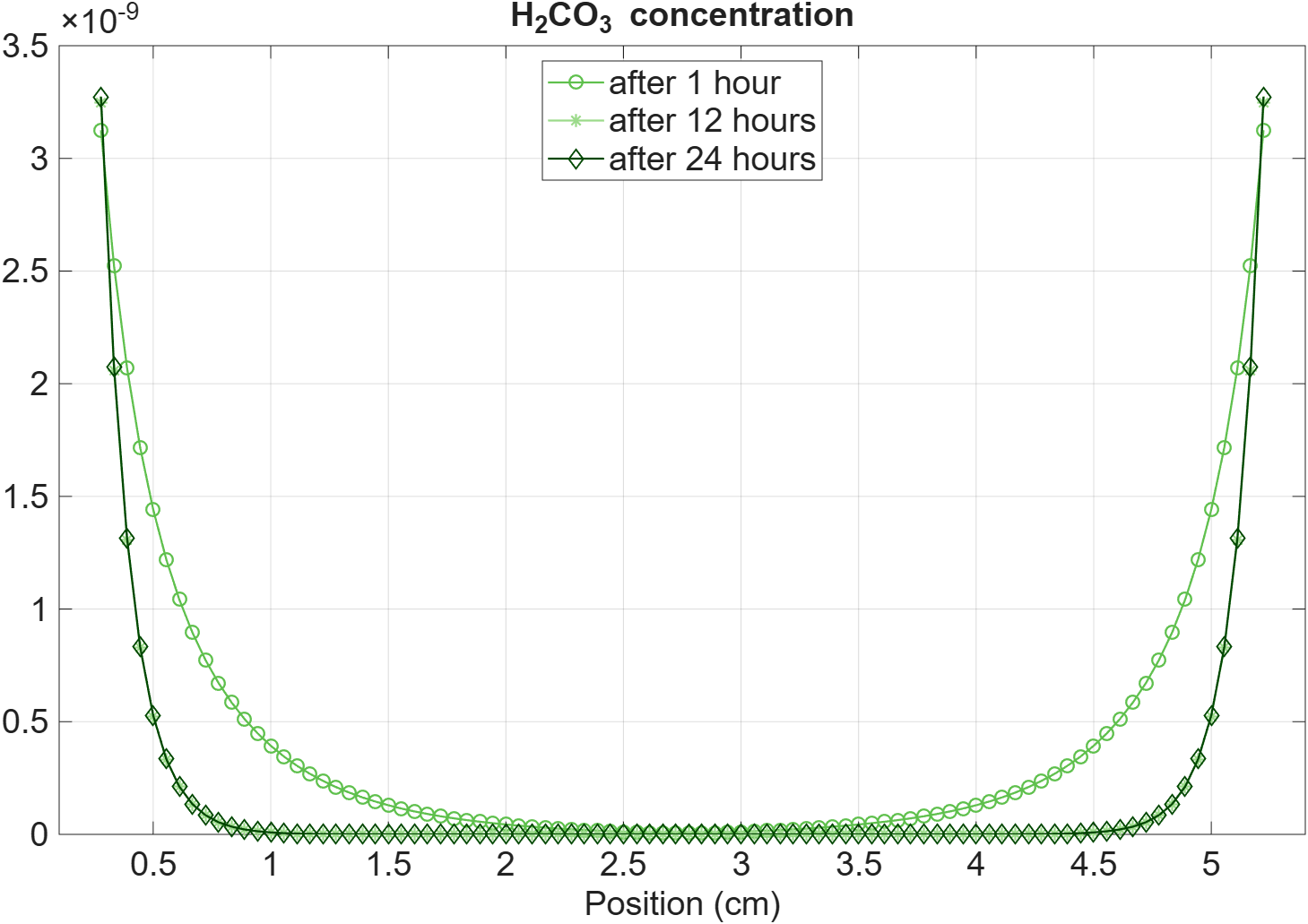}\\
	\vspace{1mm}
	\includegraphics[width=0.49\linewidth]{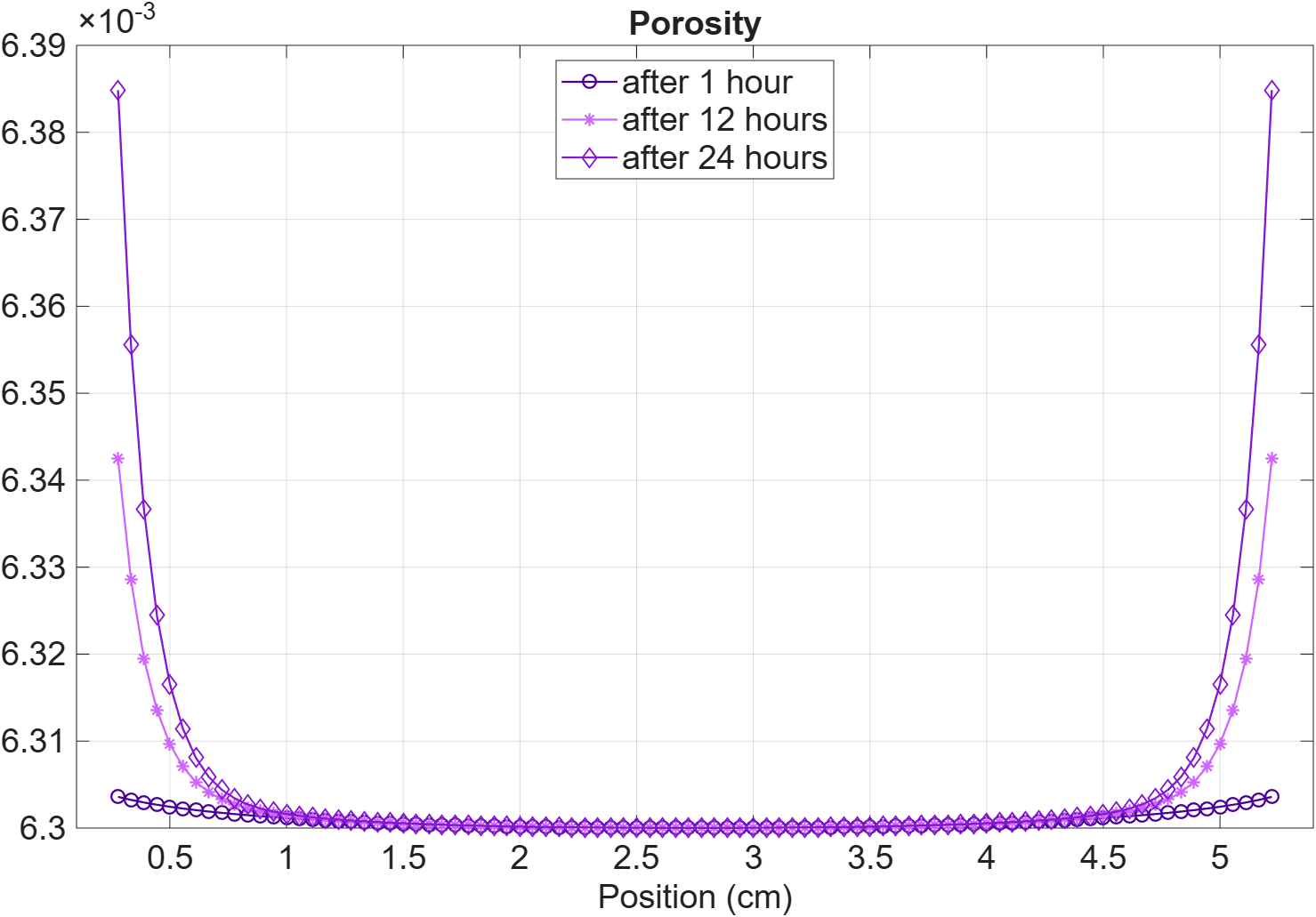}
	\includegraphics[width=0.49\linewidth]{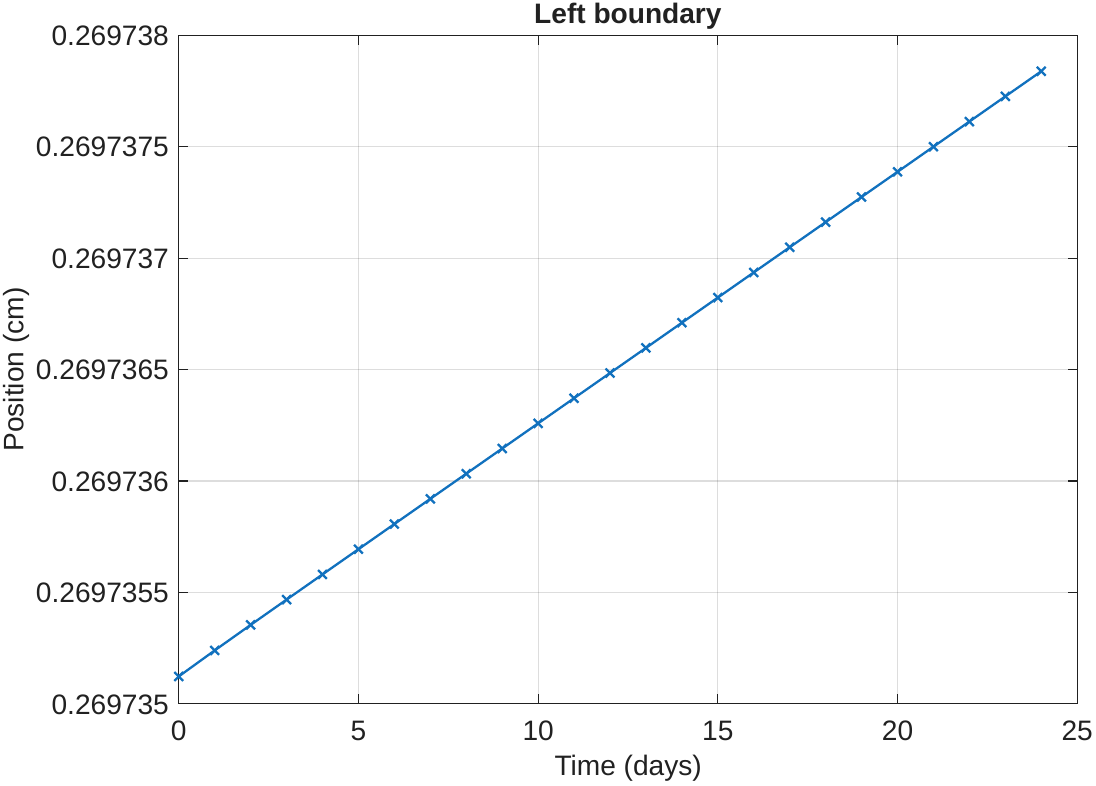}
	\caption{Results of the ``catastrophic'' case over time. First line from left: humidity and $H_2 CO_3$ concentration. Second line from left: porosity and boundary position.}
	\label{fig:catastrophe}
\end{figure}

\begin{table}
	\centering
	\resizebox{\columnwidth}{!}{%
	\begin{tabular}{|c|c|c|}
		\toprule
		Time (h) & Erosion of the left edge [cm] & Erosion of the right edge [cm]\\
		\midrule
		6 & 0.0198011446561197 & 0.0198011446561202 \\
		12 & 0.0198018272171412 & 0.0198018272171412\\
		18 & 0.0198025092994422 & 0.0198025092994421\\
		24 & 0.0198031911983033 & 0.0198031911983039\\
		\midrule
		Recession rate & 2.2919$\cdot 10^{-7}$ $cm\ s^{-1}$& 2.2919$\cdot 10^{-7}$ $cm\ s^{-1}$\\
		\bottomrule
	\end{tabular}
	}
	\caption{Slope of the erosion front for the catastrophic simulation.}
	\label{table:slopeErosion}
\end{table}
As last experiment we tested the model in a ``catastrophic'' environment that simulates the penetration of the water in a submerged stone. For this case we put $\mathcal{E}(\boldsymbol{x},t)=\widetilde{n}$ to have full saturation on the boundary, leaving the rest of the parameters unchanged. For simplicity we perform a simulation over one day, with $\Delta t = 1$ second and we use the symmetric $B(s)$. Results are shown in Figure \ref{fig:catastrophe}. As expected the water penetrates very quickly the stone reaching full saturation around 1 day. However, due to the low diffusion of the $CO_2$ from exterior we see that the stone is partially filled with $H_2 CO_3$, but then it is consumed by the reaction and only along the boundaries an accumulation can be seen. In Table \ref{table:slopeErosion} we report the recession rate, which appears to be 2.2919$\cdot 10^{-7}$ $cm\ s^{-1}$ for this experiment, which is approximately 349.76 times greater than the value 6.5528$\cdot 10^{-10} cm \ s^{-1}$ from Table \ref{table:erosion} and this is due to significantly higher humidity levels to which the stone is exposed.
\subsection{Local sensitivity analysis}
\begin{figure}
	\centering
	\includegraphics[scale=0.35]{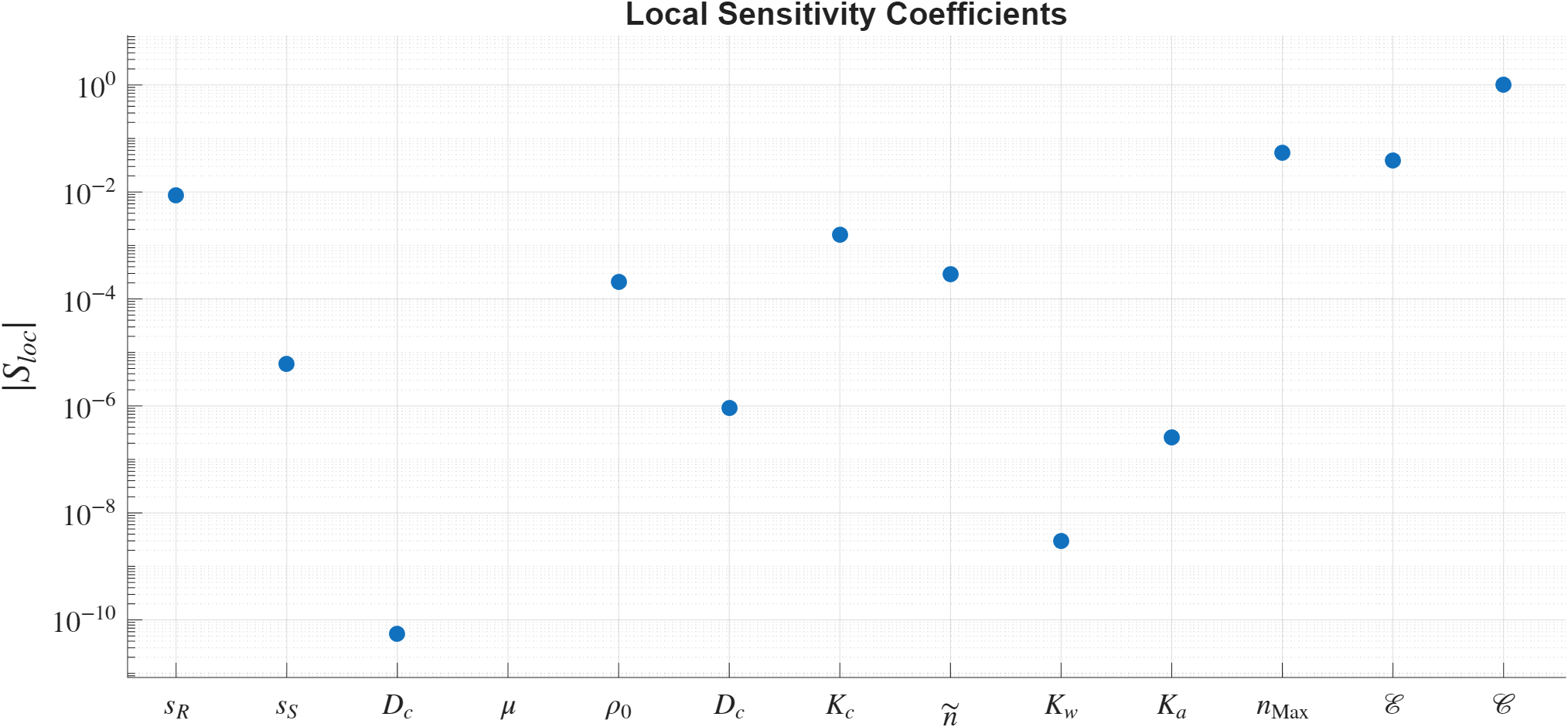}
	\caption{The absolute value of the normalized local sensitivity coefficients for recession rate function $r(t,P_{0})$ at $t=30 \ days$.}
	\label{fig:sensitivity}
\end{figure}

In order to investigate the influence of model parameters in regards to the erosion front, we conduct a local sensitivity analysis around the parameters indicated in table \ref{table:stone},\ref{table:fitting},\ref{table:kappa}.
Let $x_{B}\left(P,t\right)$ be x-coordinate of the left boundary with respect to time $t$ and parameter set $P$ and we define the recession rate function as
\begin{eqnarray}
	r(t,P):=\frac{x_{B}\left(t+\triangle t\right)-x_{B}\left(t\right)}{\triangle t}
\end{eqnarray}
with $\triangle t>0$.\\
Then we define the normalized local sensitivity coefficient for the recession rate as
\begin{eqnarray}
	S_{loc}\left(t\right):=\frac{P_{0}^{(j)}}{r(t,P_{0})}\frac{\partial r(t,P_{0})}{\partial P_{0}^{(j)}}
\end{eqnarray}
with reference parameter vector $P_{0}\in\mathbb{R}^{q}$, j-th coordinate vector $P_{0}^{(j)}:=e_{j}P_{0}$ with j-th unit vector $e_{j}$.\\
The results of the normalized local sensitivity coefficients for time $t=30\ days$ are reported in Table \ref{table:sensitivity} and depicted in Figure \ref{fig:sensitivity}.

\begin{table}
	\centering
	\begin{tabular}{|c|c|}
        \toprule
		Parameter & Local sensitivity coefficient $S_{\text{loc}}$\\
		\midrule
		$s_{R}$ & $8.723\cdot 10^{-3}$\\
		$s_{S}$ & $6.024\cdot 10^{-6}$\\
		$D_{sym}$ & $-5.5\cdot 10^{-11}$\\
		$D$ & $9.08\cdot 10^{-7}$ \\
		$\rho_{0}$ & $-2.09\cdot 10^{-4}$\\
		$\tilde{n}$& $2.88\cdot 10^{-4}$\\
		$\mu$ & 0\\
        $K_{c}$ & $1.57\cdot 10^{-3}$\\
		$K_{w}$ & $-3.017\cdot 10^{-9}$\\
		$K_{a}$ & $2.55\cdot 10^{-7}$\\
		$\mathcal{E}$ & $3.88\cdot 10^{-2}$\\
		$\mathcal{C}$ & $1.01\cdot 10^{0}$\\
		$n_{max}$ & $-5.44\cdot 10^{-2}$\\
        \bottomrule
	\end{tabular}
	\caption{Normalized local sensitivity coefficients of recession rate function $r(t,P)$ for time $t=30 \ days$ also represented in Figure \ref{fig:sensitivity}.}
	\label{table:sensitivity}
\end{table}

\begin{figure}
	\centering
	\includegraphics[width=\linewidth]{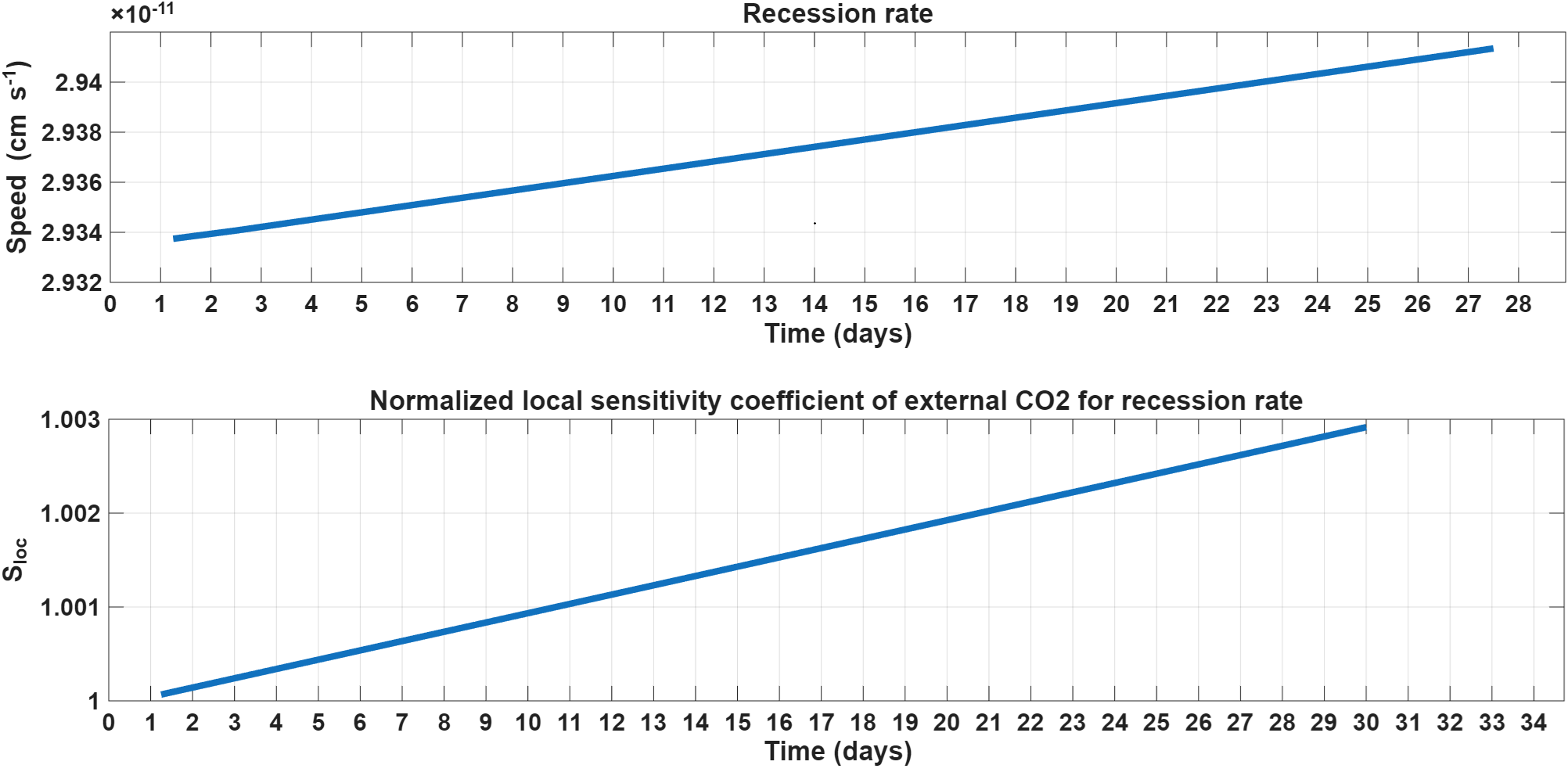}
	\caption{Normalized local sensitivity coefficients for recession rate function $r(t,P_{0})$ in time for parameter $\mathcal{C}$.}
	\label{fig:sensitivityco2}
\end{figure}

It can be seen that the external $CO_2$ concentration $\mathcal{C}$ has the largest influence on the recession rate, however $\mathcal{E}$ and $n_{max}$ do have a not negligible contribution respect to the other model parameters.\\
We need to remark that the normalized local sensitivity coefficient for the recession rate function at $P_{0}$ is for $t=30 \ days$. In Figure \ref{fig:sensitivityco2} we can observe a slight increase in the recession rate over time, accompanied by a corresponding rise in the normalized local sensitivity coefficient with respect to external $\text{CO}_2$ concentration $\mathcal{C}$. This behavior indicates that the erosion front undergoes a transient acceleration phase during the initial stages of weathering. Consequently, as the recession process accelerates, the system's dynamics become increasingly sensitive to variations in $\mathcal{C}$ over time.

\subsection{Discussion on numerical tests}

The proposed model successfully captures the dynamics of the erosion front as can be observed in the output of 1D and 2D simulations, showing an evolving interface whose speed and shape are consistent with values reported in literature \cite{brimblecombe}. 
More in detail, the computed penetration profiles of moisture and dissolved carbon dioxide qualitatively match known linear progression, see Fig. \ref{fig:asymSim}. Moreover, the advancing dissolution front over months shows realistic recession of the stone surface (see Figures \ref{fig:boundary}-\ref{fig:2Dsim}) which shows to be of several magnitude faster for the case when the stone is submerged.  We can observe that both in 1D and 2D simulations, the system reaches a steady state while the porosity of the carbonate matrix grows until it breaks the material. These results should be interpreted within the effective, dissolution-dominated regime assumed by the model: they describe the coupled advance of $CO_2$-derived aqueous acidity and porosity, but do not resolve changes in pH, carbonate speciation, calcium accumulation, saturation state or possible reprecipitation.

Regarding the order if convergence of the approximation scheme, as expected, we get a global numerical order of 1, see Table \ref{table:order}. 
Finally, we point out that the algorithm can be easily extended to 3D case.

\section{Conclusions} \label{sec:conclusions}
In this work, we have developed a mathematical framework that reproduces erosion-front advancement in carbonate stones caused by the penetration of $CO_2$-derived acidity through the water-filled pore network. Atmospheric $CO_2$ enters the model as an environmental boundary forcing, while the internal dynamics describe effective aqueous transport, reaction-driven calcite loss and the associated evolution of porosity.
The model describes the dynamics of the erosion front in accordance with literature available data, both qualitatively and quantitatively.
Importantly, the framework is flexible: different real-time environmental conditions (environmental $CO_2$ concentration, rainfall, humidity) and parameters such as temperature or material features (porosity, permeability) can easily be incorporated by setting boundary conditions and coefficients in accordance to available data. This means the model can simulate different scenarios, such as wetter vs. drier climates or higher vs. lower $CO_2$ levels, and thus can predict how these factors accelerate or slow down stone degradation for developing site-specific conservation strategies. Moreover, the level-set formulation naturally handles complex 2D geometries and moving boundaries, making it applicable to real artifacts of arbitrary shape. 
The results demonstrate that the model and numerical method provide a realistic representation of $CO_2$-driven stone erosion.
The present results should nevertheless be interpreted as a first-order macroscale representation within the effective, dissolution-dominated regime assumed by the model. The formulation does not explicitly resolve gas-liquid $CO_2$ transfer, pH, alkalinity, complete aqueous carbonate speciation, calcium and reaction-product transport, calcite saturation, or precipitation and reprecipitation. Consequently, its applicability may be limited when reaction products accumulate, the pore solution approaches calcite saturation, or precipitation becomes significant.
From the point of view of the practical application, these results should be regarded as preliminary, due to the lack of available data to calibrate the parameters $K_w$, $K_a, K_n, n_{max}$ and $\mathcal{E}$. 

This leads to a number of interesting future developments of our study, including those listed below:

\begin{itemize}
	\item Performing experiments with different carbonate stones to see the effects of the erosion for various building materials;
	\item Extending the simulation algorithm to consider 3D case;
	\item Including the gravity effect in the capillary rise; 
	\item Extending the chemical model to include explicit gas-liquid $CO_2$ transfer, carbonate acid-base equilibria, calcium and reaction-product transport, saturation-dependent dissolution and precipitation kinetics, and the possible formation of surface crusts;
	\item Assuming variable $\mathcal{E}(\boldsymbol{x},t)$ and $\mathcal{C}(\boldsymbol{x}, t)$ in space and time, aimed at simulating the effects of different climate conditions over time. Furthermore, the effect of freeze-thaw cycles  could be investigated;
	
	\item Increasing the order of the numerical scheme. Although higher order of accuracy in time can be easily achieved with a Crank-Nicolson method, a different approach should be applied in space, due to the degenerate behavior of the PDE related to the humidity diffusion.
\end{itemize}

\section*{Acknowledgements}
G.B.\, E. C. B., S. F. and M. S. are members of the Gruppo Nazionale Calcolo Scientifico-Istituto Nazionale di Alta Matematica (GNCS-INdAM).\\
The work of G.B.\, E. C. B., S. F. and M. S. has been funded by the PRIN-PNRR project MATHPROCULT Prot.\ P20228HZWR, CUP B53D23015940001.

\printbibliography
\end{document}